\documentclass[reqno]{amsart}
\usepackage{url}
\newtheorem{theorem}{Theorem}[section]

\newtheorem{corollary}[theorem]{Corollary}

\theoremstyle{definition}
\newtheorem{definition}[theorem]{Definition}

\theoremstyle{remark}

\numberwithin{equation}{section}
\usepackage{algorithmic}
\usepackage[top=3cm,bottom=2cm,right=2cm,left=2cm]{geometry}

\begin{document}

\title[An Implementation of Radu's Algorithm]{On the Computation of Identities Relating Partition Numbers in Arithmetic Progressions with Eta Quotients: An Implementation of Radu's Algorithm}

%    Remove any unused author tags.

%    author one information
\author{}
\address{}
\curraddr{}
\email{}
\thanks{}

%    author two information
\author{Nicolas Allen Smoot}
\address{}
\curraddr{}
\email{}
\thanks{}


\date{}

\dedicatory{}

\begin{abstract}
In 2015 Cristian-Silviu Radu designed an algorithm to detect identities of a class studied by Ramanujan and Kolberg, in which the generating functions of a partition function over a given set of arithmetic progression are expressed in terms of Dedekind eta quotients over a given congruence subgroup.  These identities include the famous results by Ramanujan which provide a witness to the divisibility properties of $p(5n+4),$ $p(7n+5)$.  We give an implementation of this algorithm using Mathematica.  The basic theory is first described, and an outline of the algorithm is briefly given, in order to describe the functionality and utility of our package.  We thereafter give multiple examples of applications to recent work in partition theory.  In many cases we have used our package to derive alternate proofs of various identities or congruences; in other cases we have improved previously established identities.
\end{abstract}

\maketitle

\section{Introduction}

Given some $n\in\mathbb{Z}_{\ge 0}$, we define a partition of $n$ as a weakly decreasing sequence of positive integers which sum to $n$.  Thus, the number 4 has 5 different partitions: $4,\ 3+1,\ 2+2,\ 2+1+1,\ 1+1+1+1$.  We define $p(n)$ as the number of partitions of $n$.  Thus, $p(4)=5$ (we define $p(0)=1$).

The function $p(n)$ has been seriously studied since 1748 \cite{Euler}, when Euler identified the generating function for $p(n)$ (with $q$ a formal indeterminate):

\begin{align}
\sum_{n=0}^{\infty}p(n)q^n = \prod_{m=1}^{\infty}\frac{1}{1 - q^m}.
\end{align}  However, almost nothing was known of the arithmetic properties of $p(n)$ before the twentieth century.  One of the first major breakthroughs in this area came from Ramanujan \cite{Ramanujan}:

\begin{theorem}

\begin{align}
\sum_{n=0}^{\infty}p(5n+4)q^n &= 5\cdot\prod_{m=1}^{\infty}\frac{(1-q^{5m})^5}{(1-q^m)^6},\label{Ramp5}\\
\sum_{n=0}^{\infty}p(7n+5)q^n &= 49q\cdot\prod_{m=1}^{\infty}\frac{(1-q^{7m})^7}{(1-q^m)^8} + 7\cdot\prod_{m=1}^{\infty}\frac{(1-q^{7m})^3}{(1-q^m)^4}.\label{Ramp7}
\end{align}
\end{theorem}

These are among the most iconic results in partition theory.  They are particularly interesting in that they reveal arithmetic information about $p(n)$:

\begin{theorem}
For all $n\in\mathbb{Z}_{\ge 0}$,
\begin{align*}
p(5n+4) &\equiv 0\pmod{5},\\
p(7n+5) &\equiv 0\pmod{7}.
\end{align*}
\end{theorem}  These identities are also immensely useful in the derivation of still deeper arithmetic information about $p(n)$, of which Ramanujan's classic infinite congruence families for powers of 5 and 7 serve as remarkable examples \cite{Watson}.  Moreover, the overall form of these identities conveys a deep relationship between partition numbers in arithmetic progressions and the underlying theory of modular functions.  Many interesting questions remain about this relationship.

Nearly 40 years later, Kolberg realized \cite{Kolberg} that these identities of Ramanujan could, with a very slight generalization, be extended to include a much larger variety of similar identities for $p(5n+j)$, $p(7n+j)$, $p(3n+j)$, $p(2n+j)$, and others.  For instance, Kolberg proved

\begin{align}
&\left(\sum_{n = 0}^{\infty}p(5n+1)q^n\right)\left(\sum_{n = 0}^{\infty}p(5n+2)q^n\right) = 25q\prod_{m=1}^{\infty}\frac{(1-q^{5m})^{10}}{(1-q^m)^{12}}+ 2\prod_{m=1}^{\infty}\frac{(1-q^{5m})^{4}}{(1-q^m)^{6}},
\end{align}

\begin{align}
\left(\sum_{n = 0}^{\infty}p(7n+1)q^n\right)&\left(\sum_{n = 0}^{\infty}p(7n+3)q^n\right)\left(\sum_{n = 0}^{\infty}p(7n+4)q^n\right)\\ 
=& 117649q^4\prod_{m=1}^{\infty}\frac{(1-q^{7m})^{21}}{(1-q^m)^{24}}+50421q^3\prod_{m=1}^{\infty}\frac{(1-q^{7m})^{17}}{(1-q^m)^{20}}+8232 q^2\prod_{m=1}^{\infty}\frac{(1-q^{7m})^{13}}{(1-q^m)^{16}}\nonumber\\ &+588q \prod_{m=1}^{\infty}\frac{(1-q^{7m})^{9}}{(1-q^m)^{12}}+15 \prod_{m=1}^{\infty}\frac{(1-q^{7m})^{5}}{(1-q^m)^{8}},\nonumber
\end{align}

\begin{align}
\left(\sum_{n = 0}^{\infty}p(3n)q^n\right)&\left(\sum_{n = 0}^{\infty}p(3n+1)q^n\right)\left(\sum_{n = 0}^{\infty}p(3n+2)q^n\right) \\
=& 9q\prod_{m=1}^{\infty}\frac{(1-q^{3m})(1-q^{9m})^{6}}{(1-q^m)^{10}} + 2 \prod_{m=1}^{\infty}\frac{(1-q^{3m})(1-q^{9m})^{3}}{(1-q^m)^{7}},\nonumber
\end{align}

\begin{align}
&\left(\sum_{n = 0}^{\infty}p(2n)q^n\right)\left(\sum_{n = 0}^{\infty}p(2n+1)q^n\right) = \prod_{m=1}^{\infty}\frac{(1-q^{2m})^2(1-q^{8m})^{2}}{(1-q^m)^{5}(1-q^{4m})},
\end{align} among many others.

In recent years a very large number of such identities have been produced.  They often concern the coefficients of various $q$-Pochhammer quotients, many of which can be used to enumerate various restricted partitions.

There are many different approaches by which these sorts of identities may be derived.  Kolberg, for example, proved each of the examples above (including Ramanujan's results) by manipulation of certain formal power series.  We will study them using the theory of modular functions, in the manner pioneered by Rademacher \cite{Rademacher}.

The principle behind these identities is that by changing variables to $q=e^{2\pi i\tau}$, $\tau\in\mathbb{H}$, the generating function for $p(n)$ is (very nearly) the multiplicative inverse of the Dedekind $\eta$ function.  This allows us to isolate and express the series

\begin{align*}
\sum_{n=0}^{\infty}p(mn+j)q^n,\ j,m\in\mathbb{Z}_{\ge 0},\ 0\le j < m,\ 1\le m.
\end{align*} in terms of linear combinations of $\eta$ (with a fractional input).  We can then take advantage of the symmetric properties of $\eta$ to construct a modular function using $\sum_{n\ge 0} p(mn+j)q^n$, associated with an appropriately chosen congruence subgroup $\Gamma_0(N)$.  Finally, we can ask whether this modular function is a linear combination of suitably defined eta quotients, by manipulating and studying its behavior near the boundary of $\mathbb{H}$.

What makes this a particularly powerful approach from a computational standpoint is that certain results from the theory of Riemann surfaces impose finiteness conditions on the behavior of any modular function near the boundary of $\mathbb{H}$.  This allows us to check the equality of two given modular functions by checking the equality of only a finite number of coefficients.  

Cristian-Silviu Radu recognized \cite{Radu} that this approach could be used to construct an algorithm to compute identities in the form of those discovered by Ramanujan and Kolberg above.  Indeed, he designed an algorithm which takes any arithmetic function $a(n)$ with generating function

\begin{align*}
\sum_{n=0}^{\infty}a(n)q^n = \prod_{\delta | M}\prod_{m=1}^{\infty}(1-q^{\delta m})^{r_{\delta}},
\end{align*} with $r_{\delta}\in\mathbb{Z}$ for all $\delta | M$, and an arithmetic progression $mn+j$, with $0\le j\le m-1$.  From here, and an appropriately chosen $N\in\mathbb{Z}_{\ge 2}$, the algorithm attempts to produce a set $P_{m,r}(j)\subseteq\{0, 1, 2, ..., m-1\}$ with member $j$; an integer-valued vector $s=(s_{\delta})_{\delta| N}$; and some $\alpha\in\mathbb{Z}$ such that

\begin{align}
f_{LHS} :=& f_{LHS}(s,N,M,r,m,j)(\tau)\nonumber\\
=& q^{\alpha}\prod_{\delta | N}\prod_{n=1}^{\infty}(1 - q^{\delta n})^{s_{\delta}}\cdot \prod_{j'\in P_{m,r}(j)}\sum_{n=0}^{\infty}p(mn+j')q^n
\end{align} is a modular function with certain restrictions on its behavior on the boundary of $\mathbb{H}$.

From here, we can construct a basis for the $\mathbb{Q}$-algebra generated by all eta quotients which exhibit similar behavioral restrictions to $f_{LHS}$.  We can then check membership of $f_{LHS}$ in this algebra through a finite computational means.

This paper summarizes our successful implementation of Radu's algorithm. Section 2.1 will provide a very brief review of the basic theory, and in Sections 2.2 to 2.3 an outline of our software package's structure will be given, following the design of Radu's algorithm.  Due to matters of space, we cannot provide more than a short description of the algorithm, or the underlying theory.  We highly recommend that this paper be read as a companion to \cite{Radu} and \cite{Radu0}.  We have changed the notation of these papers: notably, we denote an arithmetic progression with the letters $m,j$, rather than $m,t$ to avoid confusion with the use of $t$ as a modular function.  We have also denoted by $h_{m,j}$ what would be referred to in \cite{Radu} and \cite{Radu0} as $g_{m,j}$, and have generally reserved the letter $g$ to denote an arbitrary eta quotient.

In addition to some small notational changes from Radu's original work, we have also designed separate procedures, which account for various theoretical or computational difficulties.  We discuss these separate procedures in Section 2.4.

Apart from a description of the basic features of our package, the bulk of our paper will be examples computed by our software package.  We cover the classic cases of Ramanujan, Kolberg, and Zuckerman in Sections 3.1, 3.2.  In Section 3.3-3.4 we show examples which Radu has previously computed, and which we have given slight improvements to.  In Sections 3.5-3.9 we give applications of our package to recently discovered identities and congruences.  In many cases we are able to improve previous results.

All of the examples in this paper can be found in the Mathematica supplements available at \url{https://www3.risc.jku.at/people/nsmoot/RKAlg/RKSupplement1.nb} and \url{https://www3.risc.jku.at/people/nsmoot/RKAlg/RKSupplement2.nb}.  Section 4 explains the availability of the package, as well as its installation.

\section{Background}

\subsection{Basic Theory}

We denote $\mathbb{H}$ as the upper half complex plane, and we let $q=e^{2\pi i\tau}$, with $\tau\in\mathbb{H}$ (except in Section 3.10, wherein we will use $z\in\mathbb{H}$).  Hereafter, we will use the notation

\begin{align*}
(q^a;q^b)_{\infty} := \prod_{m=0}^{\infty}(1-q^{bm+a}).
\end{align*}  In particular,

\begin{align}
\sum_{n=0}^{\infty}p(n)q^n = \frac{1}{(q;q)_{\infty}}.
\end{align}  We now give a very brief preliminary for an understanding of the RK algorithm and its underlying theory.  We choose an approach from the perspective of compact Riemann surfaces.  We will see that this approach can yield powerful and remarkable insights on the construction and form of our RK computations, from topological properties of the associated Riemann surfaces, i.e., the classical modular curves.  We refer the reader to \cite{Diamond} for a more comprehensive treatment.  For a more classical approach, \cite{Lehner}, and \cite{Knopp} are recommended.

Let $\hat{\mathbb{H}}:=\mathbb{H}\cup\{\infty\}\cup\mathbb{Q}$.  We also define $\hat{\mathbb{Q}}:=\mathbb{Q}\cup\{\infty\}$, with $a/0=\infty$ for any $a\neq 0$.  We denote $\mathrm{SL}(2,\mathbb{Z})$ as the set of all $2\times 2$ integer matrices with determinant~1.

For any given $N\in\mathbb{Z}_{\ge 1}$, let

\begin{align*}
\Gamma_0(N) = \Bigg\{ \begin{pmatrix}
  a & b \\
  c & d 
 \end{pmatrix}\in \mathrm{SL}(2,\mathbb{Z}) : N|c \Bigg\}.
\end{align*}

We define a group action

\begin{align*}
\Gamma_0(N)\times\hat{\mathbb{H}}&\longrightarrow\hat{\mathbb{H}},\\
\left(\begin{pmatrix}
  a & b \\
  c & d 
 \end{pmatrix},\tau\right)&\longrightarrow \frac{a\tau+b}{c\tau+d}.
\end{align*} If $\gamma=\begin{pmatrix}
  a & b \\
  c & d 
 \end{pmatrix}$ and $\tau\in\hat{\mathbb{H}}$, then we write
 
 \begin{align*}
 \gamma\tau := \frac{a\tau+b}{c\tau+d}.
 \end{align*}  The orbits of this action are defined as

\begin{align*}
[\tau]_N := \left\{ \gamma\tau: \gamma\in\Gamma_0(N) \right\}.
\end{align*}

\begin{definition}
For $N\in\mathbb{Z}_{\ge 1}$, the classical modular curve of level $N$ is the set of all orbits of $\Gamma_0(N)$ applied to $\hat{\mathbb{H}}$:

\begin{align*}
\mathrm{X}_0(N):=\left\{ [\tau]_N : \tau\in\hat{\mathbb{H}} \right\}
\end{align*}
\end{definition}

The group action applied to $\hat{\mathbb{H}}$ can be restricted to $\hat{\mathbb{Q}}$: that is, for every $\tau\in\hat{\mathbb{Q}}$, $[\tau]_N\subseteq\hat{\mathbb{Q}}$.  There are only a finite number of such orbits \cite[Section 3.8]{Diamond}.

\begin{definition}
For any $N\in\mathbb{Z}_{\ge 1}$, the cusps of $\mathrm{X}_0(N)$ are the orbits of $\Gamma_0(N)$ applied to $\hat{\mathbb{Q}}$.
\end{definition}

For a detailed review of the Riemann surface structure of $\mathrm{X}_0(N)$, see \cite[Chapters 2,3]{Diamond}.  We briefly add that for each $N$, $\mathrm{X}_0(N)$ possesses a unique nonnegative integer $\mathfrak{g}\left(  \mathrm{X}_0(N) \right)$ called its genus.  This number may be computed using Theorem 3.1.1 of \cite[Chapter 3]{Diamond}.  For a theoretical understanding of the connection of the genus to RK identities, see \cite{Paule3}.  For an application, see Section 3.5.1 below.

\begin{definition}\label{DefnModular}
Let $f:\mathbb{H}\longrightarrow\mathbb{C}$ be holomorphic on $\mathbb{H}$.  Then $f$ is a modular function on $\Gamma_0(N)$ if the following properties are satisfied for every $\gamma=\left(\begin{smallmatrix}
  a & b \\
  c & d 
 \end{smallmatrix}\right)\in\mathrm{SL}(2,\mathbb{Z})$:

\begin{enumerate}
\item We have $$\displaystyle{f\left( \gamma\tau \right) = \sum_{n=n_{\gamma}}^{\infty}\alpha_{\gamma}(n)q^{n\gcd(c^2,N)/ N}},$$  with $n_{\gamma}\in\mathbb{Z}$, and $\alpha_{\gamma}(n_{\gamma})\neq 0$.  If $n_{\gamma}\ge 0$, then $f$ is holomorphic at the cusp $[a/c]_N$.  Otherwise, $f$ has a pole of order $n_{\gamma}$, and principal part
 \begin{align}
 \sum_{n=n_{\gamma}}^{-1}\alpha_{\gamma}(n)q^{n\gcd(c^2,N)/ N},\label{princpartmod}
 \end{align} at the cusp $[a/c]_N$.
 \item If $\gamma\in\Gamma_0(N)$, we have $\displaystyle{f\left( \gamma\tau \right) = f(\tau)}.$
\end{enumerate}  We refer to $n_{\gamma}(f)$ as the order of $f$ at the cusp $[a/c]_N$.
\end{definition}

As one important case, for $\gamma=\left(\begin{smallmatrix}
  1 & 0 \\
  0 & 1 
 \end{smallmatrix}\right)$, we will use the notation $\alpha_{\gamma}(n) = \alpha_{\infty}(n)$, $n_{\gamma}=n_{\infty}$.

We now define the relevant sets of all modular functions:

\begin{definition}
Let $\mathcal{M}\left(\Gamma_0(N)\right)$ be the set of all modular functions on $\Gamma_0(N)$, and $\mathcal{M}^{a/c}\left(\Gamma_0(N)\right)\subset \mathcal{M}\left(\Gamma_0(N)\right)$ to be those modular functions on $\Gamma_0(N)$ with a pole only at the cusp $[a/c]_N$.  These are both commutative algebras with 1, and standard addition and multiplication \cite[Section 2.1]{Radu}.
\end{definition}

As an additional notational matter, for any set $\mathcal{S}\subseteq\mathcal{M}\left(\Gamma_0(N)\right)$, and any field $\mathbb{K}\subseteq\mathbb{C}$, define 

\begin{align*}
\mathcal{S}_{\mathbb{K}} := \left\{ f\in\mathcal{S} : \alpha_{\mathrm{\infty}}(n)\in\mathbb{K} \text{ for all }n\ge n_{\mathrm{\infty}}(f) \right\}.
\end{align*}  Also, for any set $\mathcal{S}$ of functions on $\mathbb{C}$, denote
\begin{align*}
\left<\mathcal{S}\right>_{\mathbb{K}}:=\left\{\sum_{u=1}^v r_u\cdot g_u: g_u\in \mathcal{S}, r_u\in\mathbb{K}\right\}.
\end{align*}

Due to its precise symmetry on $\Gamma_0(N)$, any modular function $f\in\mathcal{M}\left( \Gamma_0(N) \right)$ induces a well-defined function

\begin{align*}
\hat{f}&:\mathrm{X}_0(N)\longrightarrow\mathbb{C}\cup\{\infty\}\\
&:[\tau]_N\longrightarrow f(\tau).
\end{align*}

The notions of pole order and cusps of $f$ used in Definition \ref{DefnModular} have been constructed so as to coincide with these notions applied to $\hat{f}$ on $\mathrm{X}_0(N)$.  In particular, (\ref{princpartmod}) represents the principal part of $\hat{f}$ in local coordinates near the cusp $[a/c]_N$.  Notice that as $\tau\rightarrow i\infty$, we must have $\gamma\tau\rightarrow a/c$, and $q\rightarrow 0$.

Because $f$ is holomorphic on $\mathbb{H}$ by definition, any possible poles for $\hat{f}$ must be found for $[\tau]_N\subseteq\hat{\mathbb{Q}}$.  The number and order of these poles is of paramount importance to us.

We now give an extremely important result in the general theory of Riemann surfaces \cite[Theorem 1.37]{Miranda}:

\begin{theorem}
Let $\mathrm{X}$ be a compact Riemann surface, and let $\hat{f}:\mathrm{X}\longrightarrow\mathbb{C}$ be analytic on all of $\mathrm{X}$.  Then $\hat{f}$ must be a constant function.
\end{theorem}

As an immediate consequence we have

\begin{corollary}\label{modfundtheorem}
For a given $N\in\mathbb{Z}_{\ge 1}$, if $f\in\mathcal{M}\left(\Gamma_0(N)\right)$ is holomorphic at every cusp of $\Gamma_0(N)$, then $f$ must be a constant.
\end{corollary}

It is this theorem which makes the computation of RK identities a problem of computer algebra.  Given $f,g\in\mathcal{M}^{\infty}\left(\Gamma_0(N)\right)$, we can determine whether $f=g$ by comparing the principal parts of each function.  If these parts match, then $f-g$ is holomorphic at every cusp, and must be a constant.  If the constants of $f$ and $g$ also match, then $f-g=0$.

With this important background, we will now define the eta function and describe its extraordinary properties:

\begin{definition}

For $\tau\in\mathbb{H}$, let
\begin{align*}
\eta(\tau) = q^{1/24} (q;q)_{\infty} = e^{\pi i\tau/12}\prod_{n=1}^{\infty} (1 - e^{2\pi i n\tau}).
\end{align*}
\end{definition}

The centrality of the eta function to the theory of partitions is clear in that $1/\eta(\tau)$ is effectively the generating function for $p(n)$.  Moreover, the right-hand sides of (\ref{Ramp5}), (\ref{Ramp7}), and the other Ramanujan--Kolberg identities shown in the introduction suggest that the eta function must play an important role in our algorithmic procedures.  A remarkable property of $\eta(\tau)$ is that it possesses modular symmetry, as the following theorem demonstrates:

\begin{theorem}
For any $\begin{pmatrix}
  a & b \\
  c & d 
 \end{pmatrix}\in \mathrm{SL}(2,\mathbb{Z})$, we have
\begin{align*}
\eta\left( \frac{a\tau+b}{c\tau+d} \right) = \epsilon(a,b,c,d) \left( -i (c\tau+d) \right)^{1/2}\eta(\tau),
\end{align*} with $z^{1/2}$ defined in terms of its principal branch, and $\epsilon(a,b,c,d)$ a certain 24th root of unity.
\end{theorem}  The symmetry of $\eta$ enables us to construct a very large number of modular functions on $\Gamma_0(N)$.  For example, it can be shown that

\begin{align*}
\left(\frac{\eta(5\tau)}{\eta(\tau)}\right)^6\in\mathcal{M}\left(\Gamma_0(5)\right).
\end{align*}

\begin{definition}
An eta quotient on $\Gamma_0(N)$ is an object of the form

\begin{align*}
\prod_{\lambda | N}\eta(\lambda\tau)^{s_{\lambda}}\in\mathcal{M}\left(\Gamma_0(N)\right).
\end{align*}  Denote $\mathcal{E}(N)$ as the set of all eta quotients on $\Gamma_0(N)$.  We denote $\mathcal{E}^{\infty}(N) = \mathcal{E}(N)\cap \mathcal{M}^{\infty}\left(\Gamma_0(N)\right)$.
\end{definition}

It is easy to see that $\left<\mathcal{E}^{\infty}(N)\right>_{\mathbb{K}}$ fulfills the conditions of a $\mathbb{K}$-algebra.

We will want to determine whether a given $f\in\mathcal{M}\left(\Gamma_0(N)\right)$ can be expressed as a linear combination of eta quotients, i.e., whether $f\in\left<\mathcal{E}(N)\right>_{\mathbb{Q}}$.  To do this directly, we would be forced to have a complete set of generators for $\left<\mathcal{E}(N)\right>_{\mathbb{Q}}$, and to study the behavior of $f$ at each cusp of $\Gamma_0(N)$.

To simplify the problem, we introduce the following theorem:

\begin{theorem}\label{muf}
For every $N\in\mathbb{Z}_{\ge 2}$, there exists a function $\mu\in\mathcal{E}^{\infty}(N)$ which is holomorphic at every cusp of $\Gamma_0(N)$ except $\infty$.
\end{theorem}  A proof can be found in \cite[Lemma 20]{Radu}.

We want to take some $\mu\in\mathcal{E}^{\infty}(N)$ as in Theorem \ref{muf} with a minimal value for $|n_{\infty}(\mu)|$.  Then given some $f\in\mathcal{M}\left(\Gamma_0(N)\right)$, there exists a $k_1\in\mathbb{Z}_{\ge 0}$, such that $\mu^{k_1}\cdot f\in\mathcal{M}^{\infty}\left(\Gamma_0(N)\right)$.  In this case, we need only examine the single principal part of $\mu^{k_1}\cdot f$.  We take a minimal value of $|n_{\infty}(\mu)|$ so as to keep $|n_{\infty}(\mu^{k_1})|$ small.

The question of checking whether $f\in\left<\mathcal{E}(N)\right>_{\mathbb{Q}}$ is therefore equivalent to checking whether

\begin{align*}
\mu^{k_1}\cdot f\in\mathcal{M}^{\infty}\left(\Gamma_0(N)\right)_{\mathbb{Q}}\cap\left<\mathcal{E}(N)\right>_{\mathbb{Q}}.
\end{align*}  Ideally, this membership question would equate to checking whether

\begin{align*}
\mu^{k_1}\cdot f\in\left<\mathcal{E}^{\infty}(N)\right>_{\mathbb{Q}}.
\end{align*}  However, one theoretical problem persists.  We know that

\begin{align*}
\mathcal{M}^{\infty}\left(\Gamma_0(N)\right)_{\mathbb{Q}}\cap\left<\mathcal{E}(N)\right>_{\mathbb{Q}} \supseteq \left<\mathcal{E}^{\infty}(N)\right>_{\mathbb{Q}},
\end{align*} but we have not yet established that 

\begin{align*}
\mathcal{M}^{\infty}\left(\Gamma_0(N)\right)_{\mathbb{Q}}\cap\left<\mathcal{E}(N)\right>_{\mathbb{Q}} \subseteq \left<\mathcal{E}^{\infty}(N)\right>_{\mathbb{Q}}.
\end{align*}  Current evidence suggests that the two sets are equal, and we strongly suspect that this is true.  Unfortunately, we are as of yet unable to prove it.  However, Radu was able \cite[Lemma 28]{Radu} to establish a weaker theorem:

\begin{theorem}\label{muf2}
Given some $N\in\mathbb{Z}_{\ge 2}$ and a $\mu\in\mathcal{E}^{\infty}(N)$ which is holomorphic at every cusp except $\infty$, there exists some $k_0\in\mathbb{Z}_{\ge 0}$ such that
\begin{align*}
\mu^{k_0}\cdot\left(\mathcal{M}^{\infty}\left(\Gamma_0(N)\right)_{\mathbb{Q}}\cap\left<\mathcal{E}(N)\right>_{\mathbb{Q}}\right) \subseteq \left<\mathcal{E}^{\infty}(N)\right>_{\mathbb{Q}}.
\end{align*}
\end{theorem}

The ambiguity of whether $k_0=0$ will become important later.  But what is important for the time being is that an upper bound for $k_0$ is at least computable \cite[Proof of Lemma 28]{Radu}.  With the previous two theorems, in order to check whether $f\in\left<\mathcal{E}(N)\right>_{\mathbb{Q}}$, we need to calculate a $\mu\in\mathcal{E}^{\infty}(N)$ which satisfies the conditions of Theorem \ref{muf}, with minimal order at $\infty$; we can then compute $k_0, k_1$ and check whether 

\begin{align*}
\mu^{k_0+k_1}\cdot f\in\left<\mathcal{E}^{\infty}(N)\right>_{\mathbb{Q}}.
\end{align*}

\subsection{Membership Algorithm}

Suppose that for some $N\in\mathbb{Z}_{\ge 2}$ we have a function $f\in\mathcal{M}^{\infty}\left(\Gamma_0(N)\right)_{\mathbb{Q}}$.  We know that we can expand $f$ as the following:

\begin{align}
f = \frac{c(-m_1)}{q^{m_1}} + \frac{c(-m_1+1)}{q^{m_1-1}} + ... + \frac{c(-1)}{q} + c(0) + \sum_{n=1}^{\infty} c(n)q^n,
\end{align} with $c(-m_1)\neq 0$.  Here we define $\mathrm{pord}(f):=m_1$, and $\mathrm{LC}(f) := c(-m_1)$.

We now need to define an algorithm to check the potential membership of a given $f$ in $\left<\mathcal{E}^{\infty}(N)\right>_{\mathbb{Q}}.$  We can take advantage of the very precise algebra basis which $\left<\mathcal{E}^{\infty}(N)\right>_{\mathbb{Q}}$ admits.

\begin{theorem}\label{abtheorem}
For any $N\in\mathbb{Z}_{\ge 2}$, $\mathcal{E}^{\infty}(N)$ is a finitely generated monoid.  Moreover, there exist functions $t, g_1, g_2, ..., g_v\in\mathcal{M}^{\infty}\left(\Gamma_0(N)\right)$ such that

\begin{align}
&\mathrm{pord}(t) = v,\label{ab1} \\
&\mathrm{pord}(g_i) < \mathrm{pord}(g_j),\text{ for $1\le i<j\le v$}\label{ab2} \\
&\mathrm{pord}(g_i) \not\equiv \mathrm{pord}(g_j)\pmod{\mathrm{pord}(t)},\text{ for $1\le i<j\le v$}\label{ab3} \\
&\mathrm{pord}(g_i) \not\equiv 0\pmod{\mathrm{pord}(t)},\text{ for $1\le i \le v$}\label{ab4} \\
&\left< \mathcal{E}^{\infty}(N) \right>_{\mathbb{Q}} = \left< 1, g_1, ..., g_v \right>_{\mathbb{Q}[t]}\label{ab5},
\end{align}
\end{theorem}  The proof can be found in \cite[Sections 2.1, 2.2]{Radu}.  Given any $N\in\mathbb{Z}_{\ge 2}$, the corresponding monoid generators of $\mathcal{E}^{\infty}(N)$ can be computed through a terminating algorithm \cite[Lemma 25]{Radu}.\\

\noindent \underline{PROCEDURE}: \texttt{etaGenerators} (Eta Monoid Generators)\\

\noindent \underline{INPUT}:

\noindent $N\in\mathbb{Z}_{\ge 2}$\\

\noindent \underline{OUTPUT}:

\noindent $\{\mathcal{E}_1, \mathcal{E}_2, ..., \mathcal{E}_r\}$ such that $\left\{\mathcal{E}^{k_1}_1\cdot \mathcal{E}^{k_2}_2\cdot ...\cdot \mathcal{E}^{k_r}_r: k_1, k_2, ..., k_r\in\mathbb{Z}_{\ge 0}\right\} = \mathcal{E}^{\infty}(N)$.\\ \\

Similarly, the corresponding basis elements of $\left< \mathcal{E}^{\infty}(N) \right>_{\mathbb{Q}}$ can be computed through a terminating algorithm \cite[Theorem 16]{Radu}.\\

\noindent \underline{PROCEDURE}: \texttt{AB} (Algebra Basis)\\

\noindent \underline{INPUT}:

\noindent A set of modular functions $\mathcal{F} = \{f_1, f_2, ..., f_r\}\subseteq\mathcal{M}^{\infty}\left(\Gamma_0(N)\right)$\\

\noindent \underline{OUTPUT}:

\noindent $t, g_1, g_2, ..., g_v\in\mathcal{M}^{\infty}\left(\Gamma_0(N)\right)$ such that conditions (\ref{ab1})-(\ref{ab4}) are satisfied, and $\left < 1, g_1, g_2, ..., g_v \right >_{\mathbb{Q}[t]}=\left<\mathcal{F}_{\mathrm{M}}\right>_{\mathbb{Q}}$, with $\mathcal{F}_{\mathrm{M}}:=\left\{f^{k_1}_1\cdot f^{k_2}_2\cdot ...\cdot f^{k_r}_r: k_1, k_2, ..., k_r\in\mathbb{Z}_{\ge 0}\right\}$.

If we give the output of \texttt{etaGenerators[}$N$\texttt{]} as input for \texttt{AB}, then we will produce an algebra basis satisfying the criteria of Theorem \ref{abtheorem}.

We now suppose that $f = \mu^{k_0+k_1}\cdot f_0\in\mathcal{M}^{\infty}\left(\Gamma_0(N)\right)_{\mathbb{Q}}$, with $\mu$, $k_0$, and $k_1$ defined as in Theorems \ref{muf} and \ref{muf2}.

To determine whether $f_0\in\left<\mathcal{E}(N)\right>_{\mathbb{Q}}$, we need only determine whether $f\in\left<\mathcal{E}^{\infty}(N)\right>_{\mathbb{Q}}$.  We construct an algebra basis $\{t, g_1, g_2, ..., g_v \}$ of the form described in Theorem \ref{abtheorem}; we want to know whether

\begin{align}
f\in\left<\mathcal{E}^{\infty}(N)\right>_{\mathbb{Q}}= \left< 1, g_1, g_2, ..., g_v \right>_{\mathbb{Q}[t]}.\label{memberg}
\end{align}  By Corollary \ref{modfundtheorem}, we need only examine the principal parts and constant of $f$ to determine whether (\ref{memberg}) is correct.

Because the orders of the functions $g_j$ give a complete set of representatives of the residue classes modulo $v$, we know that $m_1\equiv \mathrm{pord}(g_{j_1})\pmod{v}$, for some $j_1$ with $1\le j_1\le v$.

Suppose first that $m_1\ge \mathrm{pord}(g_{j_1})$.  Let $g_{j_1}$ have the expansion

\begin{align}
g_{j_1} = \frac{b_1(-n_1)}{q^{n_1}} + \frac{b_1(-n_1+1)}{q^{n_1-1}} + ... + \frac{b_1(-1)}{q} + b_1(0) + \sum_{n=1}^{\infty} b_1(n)q^n,
\end{align} with $b_1(-n_1)\neq 0$.  Then we can write

\begin{align}
&f_1 = f - \frac{c(-m_1)}{\mathrm{LC}\left(g_{j_1}\cdot t^{\frac{m_1 - n_1}{v}}\right)}\cdot g_{j_1}\cdot t^{\frac{m_1 - n_1}{v}},\\
&\mathrm{pord}\left( f_1\right) = m_2 < m_1.
\end{align}  

Now let $m_2\equiv \mathrm{pord}(g_{j_2})\pmod{v}$, with $1\le j_2\le v$.  If $m_2\ge \mathrm{pord}(g_{j_2})$, then we may construct $f_2$ in similar fashion as to $f_1$.

In this way, we may construct a sequence of functions $(f_l)_{l\ge 1}$, with $\mathrm{pord}(f_l)\ge\mathrm{pord}(f_{l+1})$ for each $l$.  Since $\mathrm{pord}(f_l)\in\mathbb{Z}$ for each $l$, there are two possible outcomes.  The first outcome is that for some $k\in\mathbb{Z}_{>1}$ we will have $f_{k-1}$ with

\begin{align}
f_{k-1} = c_{k-1}(0) + \sum_{n=1}^{\infty} c_{k-1}(n)q^n.
\end{align}  Of course, $c_{k-1}(0)\in\left< 1, g_1, g_2, ..., g_v \right>_{\mathbb{Q}[t]}$, so that $f_{k} = f_{k-1} - c_{k-1}(0)$ has no principal part and no constant.  In this case, we have shown that the principal part of $f$ can be constructed through combinations of the principal parts of monomials within our basis.  Since we only need to match the principal parts and constants, we can conclude that $f\in\left< 1, g_1, g_2, ..., g_v \right>_{\mathbb{Q}[t]}$.

The second outcome is that for some $k\in\mathbb{Z}_{>1}$, we produce a function $f_k$ such that $0< \mathrm{pord}(f_k) = m_k < \mathrm{pord}(g_{j_k})$.  In this case, no nonnegative power of $g_{j_k}$ can reduce the order of $f_k$, and no other element in our basis can have a matching order modulo $v$.  We must immediately conclude that the principal part of $f$ cannot be reduced in terms of the principal parts of $\left< 1, g_1, g_2, ..., g_v \right>_{\mathbb{Q}[t]}$.  Of course, this implies that $f\not\in\left< 1, g_1, g_2, ..., g_v \right>_{\mathbb{Q}[t]}$.

As we reduce the principal part of $f$, we can collect the terms 

\begin{align*}
\frac{c(-m_l)}{\mathrm{LC}\left(g_{j_l}\cdot t^{\frac{m_l - n_l}{v}}\right)}\cdot g_{j_l}\cdot t^{\frac{m_l - n_l}{v}}
\end{align*} into a set $\mathcal{V}$ of $v$ polynomials, each a sum of all the terms which use the same element $g_{j_l}$.  In the event that we can completely reduce the principal part of $f$, $\mathcal{V}$ represents the basis decomposition of $f$ over $\left< 1, g_1, g_2, ..., g_v \right>_{\mathbb{Q}[t]}$.  Below, let $\mathrm{Princ}(f)$ be the principal part of $f$ (including its constant):\\

\noindent \underline{PROCEDURE}: \texttt{MW} (Membership Witness)\\

\noindent \underline{INPUT}:

\begin{itemize}
\item $N\in\mathbb{Z}_{\ge 2}$,
\item $t,g_1,g_2,...,g_v\in\mathcal{M}^{\infty}\left(\Gamma_0(N)\right)$ satisfying (\ref{ab1})-(\ref{ab5}),
\item $\mathrm{Princ}(f)$, for some $f\in\mathcal{M}^{\infty}\left(\Gamma_0(N)\right)$.
\end{itemize}

\noindent \underline{OUTPUT}:

\begin{flalign*}
&\text{IF } f\in\left< \mathcal{E}^{\infty}(N) \right>_{\mathbb{Q}},\text{ RETURN } \{p_0, p_1, ..., p_k\}\subseteq \mathbb{Q}[x] \text{ such that}&\\
&f = \sum_{k=0}^v g_k\cdot p_k(t) \text{ with $g_0 = 1$};\\
&\text{ELSE, PRINT ``NO MEMBERSHIP"}.
\end{flalign*}

\subsection{Main Procedure}

The previous two sections discussed how to determine whether $f\in\left< \mathcal{E}(N) \right>_{\mathbb{Q}}$, for some modular function $f$.  We now need to construct the modular function $f_{LHS}$ discussed in Section 1.  Let us take an arithmetic function $a(n)$ with the generating function

\begin{align}
F_r(\tau) = \sum_{n=0}^{\infty}a(n)q^n = \prod_{\delta | M} (q^{\delta};q^{\delta})_{\infty}^{r_{\delta}},
\end{align} with $r=(r_{\delta})_{\delta | M}$ an integer-valued vector.  Suppose we are interested in a possible RK identity for $a(mn+j)$, with $0\le j< m$.

In \cite[Section 2]{Radu0}, Radu demonstrates that

\begin{align*}
&q^{\frac{24j + \sum_{\delta | M}\delta\cdot r_{\delta}}{24m}}\sum_{n=0}^{\infty} a(mn+j)q^n\\=&\frac{1}{m}\sum_{\lambda = 0}^{m-1} e^{-\frac{2\pi i\kappa\lambda}{24m} \left( 24j + \sum_{\delta | M}\delta\cdot r_{\delta} \right)}\prod_{\delta | M} \eta\left(\delta\cdot \frac{\tau+\kappa\lambda}{m} \right)^{r_{\delta}}\label{etafrac1},
\end{align*} with $\kappa = \mathrm{gcd}(m^2-1,24)$.  Therefore, if we define

\begin{align}
h_{m,j}(\tau) = q^{\frac{24j + \sum_{\delta | M}\delta\cdot r_{\delta}}{24m}}\sum_{n=0}^{\infty} a(mn+j)q^n,
\end{align}  then the functional equation of $\eta$ gives $h_{m,j}(\tau)$ a modular symmetry with respect to $\Gamma_0(N)$, for a suitably chosen $N\in\mathbb{Z}_{\ge 2}$.  Notice that this modular symmetry is not yet precise; ideally, we would want $h_{m,j}(\gamma\tau) = h_{m,j}(\tau)$ for all $\gamma\in\Gamma_0(N)$.  Instead, it can be shown \cite[Theorem 2.14]{Radu0} that for some 

\begin{align*}
\begin{pmatrix}
  a & b \\
  c & d 
 \end{pmatrix}\in\Gamma_0(N) \text{ with } a>0,c>0, \text{ and } \mathrm{gcd}(a,6)=1,
\end{align*}

\begin{align*}
h_{m,j}\left(\frac{a\tau + b}{c\tau + d}\right) = \rho\cdot\left( -i\left( c\tau + d \right) \right)^{\sum_{\delta | M}r_{\delta}/2} h_{m,j'}(\tau),
\end{align*} with $\rho := \rho(a,b,c,d,M,r,m,j)$ a certain root of unity, and $j'$ an integer, which can be computed precisely, with $0\le j' < m$.

Because $0\le j' \le m-1$, we only have a finite number of $h_{m,j'}$ to manipulate.  We can therefore take the product over all possible $h_{m,j'}$ that can be derived from $h_{m,j}(\gamma\tau)$, $\gamma\in\Gamma_0(N)$.  Denote the set of all possible $j'$ produced in this manner as $P_{m,r}(j)$.  Then the group action of any $\begin{pmatrix}
  a & b \\
  c & d 
 \end{pmatrix}\in\Gamma_0(N)$ will send $\prod_{j'\in P_{m,r}(j)}h_{m,j'}(\tau)$ to itself, multiplied by 

\begin{align}
\prod_{j' \in P_{m,r}(j)}\rho\cdot\left( -i\left( c\tau + d \right) \right)^{\sum_{\delta | M}r_{\delta}/2}.\label{extrafact}
\end{align}

To cancel the unwanted factors (\ref{extrafact}), we can construct a specific 

\begin{align*}
\prod_{\delta | N}\eta(\delta\tau)^{s_{\delta}},
\end{align*} with an integer-valued vector $s=(s_{\delta})_{\delta | N}$.  This product of eta factors will produce the multiplicative inverse of the factors we wish to cancel.  The vector $s$ is chosen so that 

\begin{align*}
\prod_{\delta | N}\eta(\delta\cdot\gamma\tau)^{s_{\delta}} = \prod_{j' \in P_{m,r}(j)}\rho^{-1}\cdot\left( -i\left( c\tau + d \right) \right)^{-\sum_{\delta | M}r_{\delta}/2} \prod_{\delta | N}\eta(\delta\tau)^{s_{\delta}}.
\end{align*}

Moreover, we also adjust $s$ so as to push the order of $\prod_{j'}h_{m,j'}$ at every cusp of $\Gamma_0(N)$ to the nonnegative integers.  That is, we incorporate the function $\mu^{k_1}$ into our system $s$.  The reasoning behind this will become clear shortly.

We can obtain $s$ as the solution to a system of equations and inequalities found in \cite[Theorems 45, 47]{Radu}.  Such a vector is guaranteed to exist for an appropriately chosen $N\in\mathbb{Z}_{\ge 2}$ \cite[Lemma 48]{Radu}.

Multiplying this eta quotient by our product of $h_{m,j'}$ factors, we obtain

\begin{align}
f_{LHS}:=& f_{LHS}(s,N,M,r,m,j)(\tau)\\ =& \prod_{\delta | N}\eta(\delta\tau)^{s_{\delta}}\cdot \prod_{j'\in P_{m,r}(j)}h_{m,j'}(\tau)\in\mathcal{M}^{\infty}\left(\Gamma_0(N)\right)_{\mathbb{Q}}.
\end{align}  We compute the set of possible solutions, and then select the optimal vector such that $f_{LHS}$ will have minimal order at $\infty$.  This is why we incorporate $\mu^{k_1}$ into our $s$ vector: doing so will greatly simplify our later calculations, since a smaller total order on the left hand side of our prospective identity ensures that less computation time will be needed to determine membership of $f_{LHS}$ (we completely ignore $\mu^{k_0}$ for the time being; see Section 2.4.3).

We now define $f_1$ as our prefactor, together with the fractional powers of $q$ taken in each $h_{m,j'}$.  This gives us another way to write $f_{LHS}$:

\begin{align}
f_1(s,N,M,r,m,j) &= \prod_{\delta | N}\eta(\delta\tau)^{s_{\delta}}\cdot q^{\sum_{j'\in P_{m,r}(j)} \frac{24j' + \sum_{\delta | M}\delta\cdot r_{\delta}}{24m}},\\
f_{LHS}(s,N,M,r,m,j) &=\prod_{\delta | N}\eta(\delta\tau)^{s_{\delta}}\cdot \prod_{j'\in P_{m,r}(j)}h_{m,j'}(\tau)\\
= f_1(s,N,&M,r,m,j)\cdot \prod_{j'\in P_{m,r}(j)}\left( \sum_{n=0}^{\infty} a(mn+j')q^n\right).
\end{align}  

At last, we come to the question of how to program $f_{LHS}$ into a computer.  Because we have previously established that $f_{LHS}$ has only one pole over $\Gamma_0(N)$, we only need to examine its principal part and constant.

Notice that $f_1$ has a principal part in $q$, and $\prod_{j'\in P_{m,r}(j)}\left( \sum_{n=0}^{\infty} a(mn+j')q^n\right)$ has no principal part in $q$.  To take the full principal part and constant of $f_{LHS}$, we need only take the principal part of $f_1$, and every term of the form $a(mn+j')q^n$, with $n\le \mathrm{pord}(f_1)$.

Let us take $\mathrm{pord}(f_1) = n_1$, and write

\begin{align*}
f_1 &= \sum_{n=-n_1}^{\infty} c(n)q^n \\
&= \frac{c(-n_1)}{q^{n_1}} + \frac{c(-n_1)+1}{q^{n_1-1}} + ... + \frac{c(1)}{q} + c(0) + \sum_{n=1}^{\infty} c(n)q^n,\\
f^{(-)}_1 &:= \frac{c(-n_1)}{q^{n_1}} + \frac{c(-n_1)+1}{q^{n_1-1}} + ... + \frac{c(1)}{q} + c(0).
\end{align*}

How many terms $a(n)$ of $F_r(\tau)$ do we need?  We know that $0\le j'\le m-1$; we also know that if $n_0> \mathrm{pord}(f_1)$, then $a(mn_0+j')q^{n_0}$ cannot contribute to the principal part of $f_{LHS}$.  Therefore, to have the principal part completely calculated, we need only take $mn+j'\le m\cdot\mathrm{pord}(f_1)+j' < m\cdot\mathrm{pord}(f_1)+m$.  We can compute and store

\begin{align*}
L &:= \sum_{n=0}^{m\cdot(\mathrm{pord}(f_1)+1)} a(n)q^n.\\
\end{align*}  We need not consider any larger values of $a(n)$.

Now define

\begin{align*}
f^{(-)}_{LHS}:= \mathrm{Princ}\left( f^{(-)}_1\cdot\prod_{j'\in P_{m,r}(j)}\left(\sum_{n=0}^{m} [mn+j']L\cdot q^n\right)\right),
\end{align*}  where by $[k]f(q)$ we mean the coefficient of $q^k$ in the expansion of $f(q)$ about $q=0$, and by $\mathrm{Princ}(f)$ we mean the principal part of $f$ together with its constant term.  We see that $f^{(-)}_{LHS}$ is a polynomial in $q^{-1}$.  In particular, $f^{(-)}_{LHS}$ is finite, and can therefore be examined by a computer.

We can now define our main procedure.  We want to determine whether our constructed $f_{LHS}\in\left< \mathcal{E}^{\infty}(N) \right>_{\mathbb{Q}}$.  We construct \cite[Section 2.1]{Radu} the functions $t,g_1, g_2, ..., g_v\in\mathcal{M}^{\infty}\left(\Gamma_0(N)\right)$, satisfying conditions (\ref{ab1})-(\ref{ab5}):

\begin{align*}
\left<\mathcal{E}^{\infty}(N)\right>_{\mathbb{Q}} = \left< 1, g_1, g_2, ..., g_v \right>_{\mathbb{Q}[t]}.
\end{align*}  We may now use our \texttt{MW} procedure to check whether $f_{LHS}\in\left< 1, g_1, g_2, ..., g_v \right>_{\mathbb{Q}[t]}$ by examining $f^{(-)}_{LHS}$.

Notice that we cannot merely construct the principal parts of the functions $t,g_l$, and disregard the rest of each function.  We reduce $f^{(-)}_{LHS}$ by subtracting monomials of the form $g_l\cdot t^n$; terms other than the principal parts of $t,g_l$ will influence the overall principal part of the product.  We must therefore be careful to construct the complete principal part of each $g_l\cdot t^n$.

If \texttt{MW} returns ``NO MEMBERSHIP", then the suspected identity does not exist---at least over $\Gamma_0(N)$.  One may attempt a different $N$ to find an identity.  Otherwise, \texttt{MW} will return

\begin{align}
\{p_0, p_1, ..., p_v\}\subseteq\mathbb{Q}[x],
\end{align} and we have the complete identity

\begin{align}
f_1(s,N,M,r,m,j)\cdot \prod_{j'\in P_{m,r}(j)}\left( \sum_{n=0}^{\infty} a(mn+j')q^n\right) = \sum_{k=0}^v g_k\cdot p_k(t).
\end{align}  Finally, we make note of an application so ubiquitous that we include it in our main procedure.  We will attempt to extract the GCD of all of the coefficients of the $p_k$.  Mathematica has a GCD procedure.  If all of the coefficients of the $p_k$ are integers, the procedure returns the GCD, which we will denote here as $\mathcal{D}$.  On the other hand, if there exists some $K\in\mathbb{Z}_{\ge 2}$ such that the coefficients are elements in $\frac{1}{K}\mathbb{Z}$, then the GCD procedure will return $\frac{1}{K}\mathcal{D}$, with $\mathcal{D}$ defined as the GCD of the coefficients with the factor $1/K$ removed.

Our procedure, \texttt{RK[}$N, M, r, m, j$\texttt{]}, takes as input an $N\in\mathbb{Z}_{\ge 2}$ which defines the congruence subgroup $\Gamma_0(N)$ to work over; a generating function (defined by $M$ and $r$), an arithmetic progression $mn+j$, with $0\le j\le m$.\\

\noindent \underline{PROCEDURE}: \texttt{RK} (Ramanujan--Kolberg Implementation)\\

\noindent \underline{INPUT}:

\begin{align}
&N\in\mathbb{Z}_{\ge 2},\label{algin4},\\
&M\in\mathbb{Z}_{\ge 1},\label{algin1}\\
&r=(r_{\delta})_{\delta | M}, r_{\delta}\in\mathbb{Z}\label{algin2}\\
&m,j\in\mathbb{Z} \text{ such that } 0\le j < m.\label{algin3}
\end{align}

\noindent \underline{OUTPUT}:

\begin{align}
&\tt{\prod_{\delta | M} (q^{\delta};q^{\delta})^{r_{\delta}}_{\infty} = \sum_{n=0}^{\infty}a(n)q^n}\label{algline1}\\
&\boxed{\tt{f_1(q)\cdot\prod_{j'\in P_{m,r}(j)}\sum_{n=0}^{\infty}a(mn+j')q^n = \sum_{g\in\mathrm{AB}} g\cdot p_g(t)}}\label{algline2}\\
&\texttt{Modular Curve: }\texttt{X}_{\texttt{0}}\texttt{(N)}\label{algline3}\\
&\texttt{N:}\ \ \ \ \ \ \ \ \ \ \ \ \ \ \ \ \ \ \ \ \ \ \ \ \ \ \ \ \ \ N\label{algline3a}\\
&\texttt{\{M, (}\texttt{r}_{\tt{\delta | M}}\texttt{)\}}\texttt{:}\ \ \ \ \ \ \ \ \ \ \ \  \{M, r\}\label{algline3b} \\
&\texttt{m :} \ \ \ \ \ \ \ \ \ \ \ \ \ \ \ \ \ \ \ \ \ \ \ \ \ \ \ \ m\label{algline3c}\\
&\texttt{P}_{\texttt{m,r}}\texttt{(j): } \ \ \ \ \ \ \ \ \ \ \ \ \ \ \ \  P_{m,r}(j)\label{algline4}\\
&\texttt{f}_{\texttt{1}}\texttt{(q): } \ \ \ \ \ \ \ \ \ \ \ \ \ \ \ \ \ \ \  f_1(q)\label{algline5}\\
&\texttt{t: }  \ \ \ \ \ \ \ \ \ \ \ \ \ \ \ \ \ \ \ \ \ \ \ \ \ t\label{algline6}\\
&\texttt{AB: } \ \ \ \ \ \ \ \ \ \ \ \ \ \ \ \ \ \ \ \ \ \ \ \ \{1, g_1, g_2, ..., g_v\}\label{algline7}\\
&\{ \texttt{p}_{\texttt{g}}\texttt{(t):g}\in\texttt{AB}\}\texttt{: }  \ \ \ \ \ \{p_{1}, p_{g_1}, ..., p_{g_v}\}\label{algline8}\\
&\texttt{Common Factor: } \ \ \ \ \  \mathcal{D}\label{algline9}
\end{align}\\

Lines (\ref{algline1}), (\ref{algline2}), (\ref{algline3}) are unsubstituted expressions which are printed before the remaining lines are computed.  They are meant to serve as a guide for the remainder of the output.  Lines (\ref{algline1}), (\ref{algline2}) indicate the form of a potential RK identity, while line (\ref{algline3}) indicates the associated modular curve.

The remaining lines give the appropriate substitutions.  First, (\ref{algline3a}), (\ref{algline3b}), (\ref{algline3c}) return $N, M, r, m$.  Line (\ref{algline4}) gives all the possible values for $j'$, including the initial input $j$.  If a vector $s$ cannot be found, then line (\ref{algline5}) will return

\begin{align*}
&\texttt{f}_{\texttt{1}}\texttt{(q): } \texttt{Select Another N}
\end{align*} indicating that we are unable to construct the necessary modular function on the given $\Gamma_0(N)$.  Similarly, if $f_{LHS}~\not\in~\left<\mathcal{E}^{\infty}(N)\right>_{\mathbb{Q}}$, then line (\ref{algline8}) will return

\begin{align*}
&\{ \texttt{p}_{\texttt{g}}\texttt{(t):g}\in\texttt{AB}\}\texttt{: } \texttt{No Membership}
\end{align*}  Otherwise, the corresponding membership witness is returned.

Finally, if a greatest common factor exists and is greater than one, then $\mathcal{D}$ is returned in line (\ref{algline9}); otherwise, the line will return

\begin{align*}
&\texttt{Common Factor: } \texttt{None}
\end{align*}

\subsection{Some Remarks}

\subsubsection{Delta}

Reexamining the identities of the introduction---(\ref{Ramp5}) and (\ref{Ramp7}) in particular---one may naturally guess that for a given progression $mn+j$, we must work over $\left<\mathcal{E}^{\infty}(m)\right>_{\mathbb{Q}}$, i.e., that the level $N$ of the associated modular curve must be equal to $m$.  In fact, while $N$ and $m$ are not always equal, they are usually closely related.  As we shall see, determination of the correct value of $N$ is an important problem for the computation of RK identities.

With a single exception, all of the exmaples found so far rely upon what Radu has called the $\Delta^{\ast}$ criterion.  For a complete definition of this criterion, see \cite[Definitions 34, 35]{Radu}.  We provide a procedure to check this criterion, in \texttt{Delta[N, M, r, m, j]}.\\

\noindent \underline{PROCEDURE}: \texttt{Delta}\\

\noindent \underline{INPUT}:  \\

\begin{align}
&N\in\mathbb{Z}_{\ge 2}\\
&M\in\mathbb{Z}_{\ge 1}\label{DEL1}\\
&r=(r_{\delta})_{\delta | M}, r_{\delta}\in\mathbb{Z}\label{DEL2}\\
&m,j\in\mathbb{Z} \text{ such that } 0\le j < m.\label{DEL3}
\end{align}  

\noindent \underline{OUTPUT}:

\begin{flalign*}
&\text{IF }\Delta^{\ast} \text{ IS SATISFIED, RETURN TRUE,}\\
&\text{ELSE, RETURN FALSE} 
\end{flalign*}\\

We provide an additional procedure, in \texttt{minN[M, r, m, j]}, which will compute the minimal $N$ that satisfies the $\Delta^{\ast}$ criterion.

\noindent \underline{PROCEDURE}: \texttt{minN}\\

\noindent \underline{INPUT}:  \\

\begin{align*}
&M\in\mathbb{Z}_{\ge 1}\\
&r=(r_{\delta})_{\delta | M}, r_{\delta}\in\mathbb{Z}\\
&m,j\in\mathbb{Z} \text{ such that } 0\le j < m.
\end{align*}  

\noindent \underline{OUTPUT}:

\begin{flalign*}
&N\in\mathbb{Z}_{\ge 2} \text{ such that } \texttt{Delta[N, M, r, m, j]=True}.
\end{flalign*}\\

The \texttt{RK} algorithm works over two distinct cases: Case 1, in which the $\Delta^{\ast}$ criterion is satisfied, and Case 2, in which it fails \cite[Section 3.1]{Radu}.  The great majority of identities we have found arise from the first case.  We will provide one interesting example of an identity arising from Case 2.  However, Case 1 is generally a faster algorithm, and we recommend that users compute an $N$ for which the $\Delta^{\ast}$ criterion is satisfied.

At any rate, for any given $M, r=(r_{\delta})_{\delta | M}, m, j$ with $0\le j < m$, there must exist an $N\in\mathbb{Z}_{\ge 2}$ such that the $\Delta^{\ast}$ criterion is satisfied \cite[Section 3.1]{Radu}.  It is generally convenient to work with the smallest possible $N$ that satisfies the criterion.  However, we will see in subsequent examples that the smallest possible case is not always the most useful.  We will therefore leave the criterion for establishing $N$ as separate from the main algorithm, and define $N$ as part of the input.

\subsubsection{\texttt{RKMan}}

We also include a slightly modified implementation that we refer to as \texttt{RKMan}.  This procedure is nearly identical to that used for Radu's algorithm, except that the algebra basis is included in the input.  This is often helpful because, as we will see in some examples, construction of the algebra basis for $\left< \mathcal{E}^{\infty}(N) \right>_{\mathbb{Q}}$ is often inefficient.  If we already have a suitable algebra basis calculated (perhaps from a database, or a general study of eta quotient spaces), and if we know the genus of the corresponding Riemann surface, we may be able to construct a basis by inspection.  This can often easily shorten the computation time.  See Section 3.5.1 for an example.

\subsubsection{\texttt{RKE}}

Regarding the value of $k_0$ in Theorem 5, we very strongly suspect that $k_0$ may always be set to 0, and that therefore 

\begin{align*}
\mathcal{M}^{\infty}\left(\Gamma_0(N)\right)_{\mathbb{Q}}\cap\left<\mathcal{E}(N)\right>_{\mathbb{Q}} = \left<\mathcal{E}^{\infty}(N)\right>_{\mathbb{Q}}.
\end{align*} for all $N\in\mathbb{Z}_{\ge 2}$.  This is important, because the computation of a bound for $k_0$ is costly, and increases the runtime of our package.  We therefore include the procedure \texttt{RKE} in addition to \texttt{RK} command.  The two commands are nearly identical, except that \texttt{RKE} includes the power $\mu^{k_0}$ in our prefactor.  This often increases runtime.  See our examples at \url{https://www3.risc.jku.at/people/nsmoot/RKAlg/RKSupplement1.nb}.

We also include the procedure \texttt{RKManE}, which is identical to \texttt{RKMan}, except that it includes $\mu^{k_0}$.

In the examples below, we use the procedures \texttt{RK}, \texttt{RKMan}.

\section{Examples}

We now give an overview of applications of our package.  Except for Sections 3.1-3.2, which cover the classic cases, each of our examples is chosen from contemporary work done in partition theory over the last ten years---in most cases, within the last five years.  In many cases we give substantial improvements on previous results, and (with the notable exception of the identities found with respect to $\bar{p}(n)$) the necessary computations take a few minutes at most on a modest laptop.

\subsection{Ramanujan's Classics}

The most obvious examples to check are the classic identities of Ramanujan and Kolberg for $p(5n+4)$ and $p(7n+5)$.

The generating function for $p(n)$ is of course $1/(q;q)_{\infty}$, which can be described by setting $M=1$, $r=(-1)$.  If we now take $m=5$, guess $N=5$, and take $j=4$, then we have

\begin{flalign*}
\texttt{In[1] = }&\texttt{RK} [ 5,1,\{-1\},5,4 ]&\\
&\tt{\prod_{\delta | M} (q^{\delta};q^{\delta})^{r_{\delta}}_{\infty} = \sum_{n=0}^{\infty}a(n)q^n}\\
&\boxed{\tt{f_1(q)\cdot\prod_{j'\in P_{m,r}(j)}\sum_{n=0}^{\infty}a(mn+j')q^n = \sum_{g\in\mathrm{AB}} g\cdot p_g(t)}}\\
&\texttt{Modular Curve: }\texttt{X}_{\texttt{0}}\texttt{(N)}\\
\texttt{Out[1] = }&\\
&\texttt{N:}\ \ \ \ \ \ \ \ \ \ \ \ \ \ \ \ \ \ \ \ \ \ \ \ \ \ \ \ \ \ 5\\
&\texttt{\{M, (}\texttt{r}_{\tt{\delta | M}}\texttt{)\}}\texttt{:}\ \ \ \ \ \ \ \ \ \ \ \  \{1, \{-1\}\} \\
&\texttt{m :} \ \ \ \ \ \ \ \ \ \ \ \ \ \ \ \ \ \ \ \ \ \ \ \ \ \ \ \ 5\\
&\texttt{P}_{\texttt{m,r}}\texttt{(j): } \ \ \ \ \ \ \ \ \ \ \ \ \ \ \ \  \{4\}\\
&\texttt{f}_{\texttt{1}}\texttt{(q): } \ \ \ \ \ \ \ \ \ \ \ \ \ \ \ \ \ \ \  \frac{((q;q)_{\infty})^6}{((q^5;q^5)_{\infty})^5}\\
&\texttt{t: }  \ \ \ \ \ \ \ \ \ \ \ \ \ \ \ \ \ \ \ \ \ \ \ \ \ \frac{((q;q)_{\infty})^6}{q((q^5;q^5)_{\infty})^6}\\
&\texttt{AB: } \ \ \ \ \ \ \ \ \ \ \ \ \ \ \ \ \ \ \ \ \ \ \ \ \{1\}\\
&\{ \texttt{p}_{\texttt{g}}\texttt{(t):g}\in\texttt{AB}\}\texttt{: }  \ \ \ \ \ \{5\}\\
&\texttt{Common Factor: } \ \ \ \ \ \ \  5
\end{flalign*}

We see that $P_{m,r}(j)=\{4\}$, indicating that our left hand side will only contain the series $\sum_{n\ge 0}p(5n+4)q^n$.  With $f_1$, we have the left hand side of any possible identity as

\begin{align*}
f_{LHS} = \frac{(q;q)_{\infty}^6}{(q^5;q^5)_{\infty}^5}\sum_{n = 0}^{\infty}p(5n+4)q^n \in\mathcal{M}^{\infty}\left(\Gamma_0(5)\right).
\end{align*}  In this case our algebra basis is extremely simple:

\begin{align*}
\left<\mathcal{E}^{\infty}(5)\right>_{\mathbb{Q}} = \left< 1 \right>_{\mathbb{Q}[t]} = \mathbb{Q}[t],
\end{align*} with

\begin{align*}
t = \frac{(q;q)_{\infty}^6}{q(q^5;q^5)_{\infty}^6}.
\end{align*}  Because the basis contains only the identity, we only need a single polynomial in $t$.  In this case, the polynomial is 5.

\begin{align*}
\frac{(q;q)_{\infty}^6}{(q^5;q^5)_{\infty}^5}\sum_{n = 0}^{\infty}p(5n+4)q^n = 5.
\end{align*}  A quick rearrangement gives us (\ref{Ramp5})

Similarly, taking $m=7$, $j=5$, and guessing $N=7$, we have

\begin{flalign*}
\texttt{In[2] = }&\texttt{RK} [7,1,\{-1\},7,5]&\\
&\tt{\prod_{\delta | M} (q^{\delta};q^{\delta})^{r_{\delta}}_{\infty} = \sum_{n=0}^{\infty}a(n)q^n}\\
&\boxed{\tt{f_1(q)\cdot\prod_{j'\in P_{m,r}(j)}\sum_{n=0}^{\infty}a(mn+j')q^n = \sum_{g\in\mathrm{AB}} g\cdot p_g(t)}}\\
&\texttt{Modular Curve: }\texttt{X}_{\texttt{0}}\texttt{(N)}\\
\texttt{Out[2] = }&\\
&\texttt{N:}\ \ \ \ \ \ \ \ \ \ \ \ \ \ \ \ \ \ \ \ \ \ \ \ \ \ \ \ \ \ 7\\
&\texttt{\{M, (}\texttt{r}_{\tt{\delta | M}}\texttt{)\}}\texttt{:}\ \ \ \ \ \ \ \ \ \ \ \  \{1, \{-1\}\} \\
&\texttt{m :} \ \ \ \ \ \ \ \ \ \ \ \ \ \ \ \ \ \ \ \ \ \ \ \ \ \ \ \ 7\\
&\texttt{P}_{\texttt{m,r}}\texttt{(j): } \ \ \ \ \ \ \ \ \ \ \ \ \ \ \ \  \{5\}\\
&\texttt{f}_{\texttt{1}}\texttt{(q): } \ \ \ \ \ \ \ \ \ \ \ \ \ \ \ \ \ \ \  \frac{((q;q)_{\infty})^8}{q((q^7;q^7)_{\infty})^7}\\
&\texttt{t: }  \ \ \ \ \ \ \ \ \ \ \ \ \ \ \ \ \ \ \ \ \ \ \ \ \ \frac{((q;q)_{\infty})^4}{q((q^7;q^7)_{\infty})^4}\\
&\texttt{AB: } \ \ \ \ \ \ \ \ \ \ \ \ \ \ \ \ \ \ \ \ \ \ \ \ \{1\}\\
&\{ \texttt{p}_{\texttt{g}}\texttt{(t):g}\in\texttt{AB}\}\texttt{: }  \ \ \ \ \ \{49+7t\}\\
&\texttt{Common Factor: } \ \ \ \ \ \ \  7
\end{flalign*}

This gives us

\begin{align*}
\frac{(q;q)_{\infty}^8}{q(q^7;q^7)_{\infty}^7}\sum_{n = 0}^{\infty}p(7n+5)q^n = 49 + 7\frac{(q;q)_{\infty}^4}{q(q^7;q^7)_{\infty}^4},
\end{align*} which yields (\ref{Ramp7}) on rearrangement.

In the following examples, we will omit the three printed lines, as well as the first three lines of output from each example for the sake of brevity.

\subsection{Classic Identities by Kolberg and Zuckerman}

A large number of classic analogues to Ramanujan's results have been found.  We start with an identity discovered by Zuckerman \cite{Zuckerman} for $p(13n+6)$.

\begin{theorem}
\begin{align*}
\sum_{n=0}^{\infty}p(13n+6)q^n =& 11\frac{(q^{13};q^{13})_{\infty}}{(q;q)_{\infty}^2} + 468q\frac{(q^{13};q^{13})_{\infty}^3}{(q;q)_{\infty}^4} + 6422q^2\frac{(q^{13};q^{13})_{\infty}^5}{(q;q)_{\infty}^6}\\& + 43940q^3\frac{(q^{13};q^{13})_{\infty}^7}{(q;q)_{\infty}^8} + 171366q^4\frac{(q^{13};q^{13})_{\infty}^{9}}{(q;q)_{\infty}^{10}}\\& + 371293q^5\frac{(q^{13};q^{13})_{\infty}^{11}}{(q;q)_{\infty}^{12}} + 371293q^6\frac{(q^{13};q^{13})_{\infty}^{13}}{(q;q)_{\infty}^{14}}.
\end{align*}
\end{theorem}

\begin{flalign*}
\texttt{In[3] = }&\texttt{RK} [13,1,\{-1\},13,6]&\\
\texttt{Out[3] = }&\\
&\texttt{P}_{\texttt{m,r}}\texttt{(j): } \ \ \ \ \ \ \ \ \ \ \ \ \ \ \ \  \{6\}\\
&\texttt{f}_{\texttt{1}}\texttt{(q): } \ \ \ \ \ \ \ \ \ \ \ \ \ \ \ \ \ \ \  \frac{((q;q)_{\infty})^{14}}{q^6((q^{13};q^{13})_{\infty})^{13}}\\
&\texttt{t: }  \ \ \ \ \ \ \ \ \ \ \ \ \ \ \ \ \ \ \ \ \ \ \ \ \ \frac{((q;q)_{\infty})^2}{q((q^{13};q^{13})_{\infty})^2}\\
&\texttt{AB: } \ \ \ \ \ \ \ \ \ \ \ \ \ \ \ \ \ \ \ \ \ \ \ \ \{1\}\\
&\{ \texttt{p}_{\texttt{g}}\texttt{(t):g}\in\texttt{AB}\}\texttt{: }  \ \ \ \ \ \{371293+371293 t+171366 t^2+43940 t^3+6422 t^4+468 t^5+11 t^6\}\\
&\texttt{Common Factor: } \ \ \ \ \ \ \  \texttt{None}
\end{flalign*}

We will now use our algorithm to derive the identities which Kolberg found \cite{Kolberg} for $p(5n+j)$, $p(7n+j)$, and $p(3n+j)$. 

Starting with $p(5n+j)$ for $0\le j\le 4$, if we take $N=5$ once more, and set $j=1$, \cite[(4.2)]{Kolberg} we have

\begin{flalign*}
\texttt{In[4] = }&\texttt{RK} [5,1,\{-1\},5,1]&\\
\texttt{Out[4] = }&\\
&\texttt{P}_{\texttt{m,r}}\texttt{(j): } \ \ \ \ \ \ \ \ \ \ \ \ \ \ \ \  \{1,2\}\\
&\texttt{f}_{\texttt{1}}\texttt{(q): } \ \ \ \ \ \ \ \ \ \ \ \ \ \ \ \ \ \ \  \frac{((q;q)_{\infty})^{12}}{((q^5;q^5)_{\infty})^{10}}\\
&\texttt{t: }  \ \ \ \ \ \ \ \ \ \ \ \ \ \ \ \ \ \ \ \ \ \ \ \ \ \frac{((q;q)_{\infty})^6}{q((q^5;q^5)_{\infty})^6}\\
&\texttt{AB: } \ \ \ \ \ \ \ \ \ \ \ \ \ \ \ \ \ \ \ \ \ \ \ \ \{1\}\\
&\{ \texttt{p}_{\texttt{g}}\texttt{(t):g}\in\texttt{AB}\}\texttt{: }  \ \ \ \ \ \{25+2t\}\\
&\texttt{Common Factor: } \ \ \ \ \ \ \  \texttt{None}
\end{flalign*}

Working over the same congruence subgroup $\Gamma_0(5)$, we keep the same algebra basis and $t$.  The most notable difference is that we have the product 

\begin{align*}
\left(\sum_{n\ge 0}p(5n+1)q^n\right)\left(\sum_{n\ge 0}p(5n+2)q^n \right)
\end{align*} on the left hand side.  Our right-hand side is given as a more complicated $25+2t$, and we have

\begin{align*}
\frac{(q;q)_{\infty}^{12}}{(q^5;q^5)_{\infty}^{10}}\left(\sum_{n = 0}^{\infty}p(5n+1)q^n\right)\left(\sum_{n = 0}^{\infty}p(5n+2)q^n\right) = 25+2\frac{(q;q)_{\infty}^6}{q(q^5;q^5)_{\infty}^6}.
\end{align*}  We can similarly examine $j=3$ \cite[(4.3)]{Kolberg} and derive the identity

\begin{align*}
\frac{(q;q)_{\infty}^{12}}{(q^5;q^5)_{\infty}^{10}}\left(\sum_{n = 0}^{\infty}p(5n+3)q^n\right)\left(\sum_{n = 0}^{\infty}p(5n)q^n\right) = 25+3\frac{(q;q)_{\infty}^6}{q(q^5;q^5)_{\infty}^6}.
\end{align*}

On the other hand, we can set $m=7, j=1, N=7$, \cite[(5.2)]{Kolberg} and we will derive

\begin{flalign*}
\texttt{In[5] = }&\texttt{RK} [7,1,\{-1\},7,1]&\\
\texttt{Out[5] = }&\\
&\texttt{P}_{\texttt{m,r}}\texttt{(j): } \ \ \ \ \ \ \ \ \ \ \ \ \ \ \ \  \{1,3,4\}\\
&\texttt{f}_{\texttt{1}}\texttt{(q): } \ \ \ \ \ \ \ \ \ \ \ \ \ \ \ \ \ \ \  \frac{((q;q)_{\infty})^{24}}{q((q^7;q^7)_{\infty})^{21}}\\
&\texttt{t: }  \ \ \ \ \ \ \ \ \ \ \ \ \ \ \ \ \ \ \ \ \ \ \ \ \ \frac{((q;q)_{\infty})^4}{q((q^7;q^7)_{\infty})^4}\\
&\texttt{AB: } \ \ \ \ \ \ \ \ \ \ \ \ \ \ \ \ \ \ \ \ \ \ \ \ \{1\}\\
&\{ \texttt{p}_{\texttt{g}}\texttt{(t):g}\in\texttt{AB}\}\texttt{: }  \ \ \ \ \ \{117649+50421 t+8232 t^2+588 t^3+15 t^4\}\\
&\texttt{Common Factor: } \ \ \ \ \ \ \  \texttt{None}
\end{flalign*}

and the identity

\begin{align*}
&\frac{(q;q)_{\infty}^{24}}{q^4(q^7;q^7)_{\infty}^{21}}\left(\sum_{n = 0}^{\infty}p(7n+1)q^n\right)\left(\sum_{n = 0}^{\infty}p(7n+3)q^n\right)\left(\sum_{n = 0}^{\infty}p(7n+4)q^n\right)\\ =& 117649+50421 \frac{(q;q)_{\infty}^4}{q(q^7;q^7)_{\infty}^4}+8232 \frac{(q;q)_{\infty}^8}{q^2(q^7;q^7)_{\infty}^8}+588 \frac{(q;q)_{\infty}^{12}}{q^3(q^7;q^7)_{\infty}^{12}}+15 \frac{(q;q)_{\infty}^{16}}{q^4(q^7;q^7)_{\infty}^{16}}.
\end{align*}  The corresponding identity for $p(7n+2)$ \cite[(5.3)]{Kolberg} can be easily found.

Finally, we set $m=3$, $j=1$, $N=9$, \cite[(3.4)]{Kolberg} and derive

\begin{flalign*}
\texttt{In[6] = }&\texttt{RK} [9,1,\{-1\},3,1]&\\
\texttt{Out[6] = }&\\
&\texttt{P}_{\texttt{m,r}}\texttt{(j): } \ \ \ \ \ \ \ \ \ \ \ \ \ \ \ \  \{0,1,2\}\\
&\texttt{f}_{\texttt{1}}\texttt{(q): } \ \ \ \ \ \ \ \ \ \ \ \ \ \ \ \ \ \ \  \frac{((q;q)_{\infty})^{10}}{q(q^3;q^3)_{\infty}((q^9;q^9)_{\infty})^{6}}\\
&\texttt{t: }  \ \ \ \ \ \ \ \ \ \ \ \ \ \ \ \ \ \ \ \ \ \ \ \ \ \frac{((q;q)_{\infty})^3}{q((q^9;q^9)_{\infty})^3}\\
&\texttt{AB: } \ \ \ \ \ \ \ \ \ \ \ \ \ \ \ \ \ \ \ \ \ \ \ \ \{1\}\\
&\{ \texttt{p}_{\texttt{g}}\texttt{(t):g}\in\texttt{AB}\}\texttt{: }  \ \ \ \ \ \{9+2t\}\\
&\texttt{Common Factor: } \ \ \ \ \ \ \  \texttt{None}
\end{flalign*}

And we have

\begin{align*}
&\frac{(q;q)^{10}}{q(q^3;q^3)(q^9;q^9)^{6}}\left(\sum_{n = 0}^{\infty}p(3n)q^n\right)\left(\sum_{n = 0}^{\infty}p(3n+1)q^n\right)\left(\sum_{n = 0}^{\infty}p(3n+2)q^n\right)\\ =& 9 + 2\frac{(q;q)^3}{q(q^9;q^9)^3}.
\end{align*}

Finally, we give another result found by Kolberg \cite[(2.4)]{Kolberg0}.  We set $m=2$, $j=1$, $N=8$ and derive

\begin{flalign*}
\texttt{In[7] = }&\texttt{RK} [8,1,\{-1\},2,1]&\\
\texttt{Out[7] = }&\\
&\texttt{P}_{\texttt{m,r}}\texttt{(j): } \ \ \ \ \ \ \ \ \ \ \ \ \ \ \ \  \{0,1\}\\
&\texttt{f}_{\texttt{1}}\texttt{(q): } \ \ \ \ \ \ \ \ \ \ \ \ \ \ \ \ \ \ \  \frac{((q;q)_{\infty})^{5}(q^4;q^4)_{\infty}}{((q^2;q^2)_{\infty})^2((q^8;q^8)_{\infty})^{2}}\\
&\texttt{t: }  \ \ \ \ \ \ \ \ \ \ \ \ \ \ \ \ \ \ \ \ \ \ \ \ \ \frac{((q^4;q^4)_{\infty})^{12}}{q((q^2;q^2)_{\infty})^4((q^8;q^8)_{\infty})^{8}}\\
&\texttt{AB: } \ \ \ \ \ \ \ \ \ \ \ \ \ \ \ \ \ \ \ \ \ \ \ \ \{1\}\\
&\{ \texttt{p}_{\texttt{g}}\texttt{(t):g}\in\texttt{AB}\}\texttt{: }  \ \ \ \ \ \{1\}\\
&\texttt{Common Factor: } \ \ \ \ \ \ \  \texttt{None}
\end{flalign*}

And we have

\begin{align*}
&\frac{((q;q)_{\infty})^{5}(q^4;q^4)_{\infty}}{((q^2;q^2)_{\infty})^2((q^8;q^8)_{\infty})^{2}}\left(\sum_{n = 0}^{\infty}p(2n)q^n\right)\left(\sum_{n = 0}^{\infty}p(2n+1)q^n\right) = 1.
\end{align*}

\subsection{Radu's Identity for 11}

A substantial amount of work has been done attempting a witness identity for $p(11n+6)\equiv 0\pmod{11}$.  We will show one interesting attempt by Radu, though we hasten to add that a great deal of work has been done by others on the problem (for an interesting approach, see \cite{Hemmecke}).  If we were to attempt to find such an identity for $M=1$, $r=(-1)$, $m=11$, $N=11$, $j=6$, then our algorithm returns

\begin{flalign*}
\texttt{In[8] = }&\texttt{RK} [11,1,\{-1\},11,6]&\\
\texttt{Out[8] = }&\\
&\texttt{P}_{\texttt{m,r}}\texttt{(j): } \ \ \ \ \ \ \ \ \ \ \ \ \ \ \ \  \{6\}\\
&\texttt{f}_{\texttt{1}}\texttt{(q): } \ \ \ \ \ \ \ \ \ \ \ \ \ \ \ \ \ \ \  \frac{(q;q)_{\infty}^{12}}{q^4(q^{11};q^{11})_{\infty}^{11}}\\
&\texttt{t: }  \ \ \ \ \ \ \ \ \ \ \ \ \ \ \ \ \ \ \ \ \ \ \ \ \ \frac{(q;q)_{\infty}^{12}}{q^5(q^{11};q^{11})_{\infty}^{12}}\\
&\texttt{AB: } \ \ \ \ \ \ \ \ \ \ \ \ \ \ \ \ \ \ \ \ \ \ \ \ \{1\}\\
&\{ \texttt{p}_{\texttt{g}}\texttt{(t):g}\in\texttt{AB}\}\texttt{: }  \ \ \ \ \ \texttt{No Membership}\\
&\texttt{Common Factor: } \ \ \ \ \ \ \  \texttt{None}
\end{flalign*}

Our membership witness returns a null result, indicating that our constructed modular function does not lie within $\left< \mathcal{E}^{\infty}(11) \right>_{\mathbb{Q}}$.

If we were to take $N=22$, however, we get

\begin{flalign*}
\texttt{In[9] = }&\texttt{RK} [22,1,\{-1\},11,6]&\\
\texttt{Out[9] = }&\\
&\texttt{P}_{\texttt{m,r}}\texttt{(j): } \ \ \ \ \ \ \ \ \ \ \ \ \  \{6\}\\
&\texttt{f}_{\texttt{1}}\texttt{(q): } \ \ \ \ \ \ \ \ \ \ \ \ \ \ \ \  \frac{(q;q)_{\infty}^{12}(q^2;q^2)_{\infty}^{2}(q^{11};q^{11})_{\infty}^{11}}{q^{14}(q^{22};q^{22})_{\infty}^{22}}\\
&\texttt{t: }  \ \ \ \ \ \ \ \ \ \ \ \ \ \ \ \ \ \ \ \ \ \  -\frac{1}{8}\frac{(q^2;q^2)_{\infty}(q^{11};q^{11})_{\infty}^{11}}{q^5(q;q)_{\infty}(q^{22};q^{22})_{\infty}^{11}}+\frac{1}{11}\frac{(q^2;q^2)_{\infty}^{8}(q^{11};q^{11})_{\infty}^{4}}{q^5(q;q)_{\infty}^{4}(q^{22};q^{22})_{\infty}^{8}} +\frac{3}{88}\frac{(q;q)_{\infty}^{7}(q^{11};q^{11})_{\infty}^{3}}{q^5(q^2;q^2)_{\infty}^{3}(q^{22};q^{22})_{\infty}^{7}}\\
&\texttt{AB: } \ \ \ \ \ \ \ \ \ \ \ \ \ \ \ \ \ \ \ \ \ \{1,-\frac{1}{8}\frac{(q^2;q^2)_{\infty}(q^{11};q^{11})_{\infty}^{11}}{q^5(q;q)_{\infty}(q^{22};q^{22})_{\infty}^{11}}+\frac{2}{11}\frac{(q^2;q^2)_{\infty}^{8}(q^{11};q^{11})_{\infty}^{4}}{q^5(q;q)_{\infty}^{4}(q^{22};q^{22})_{\infty}^{8}} +\frac{5}{88}\frac{(q;q)_{\infty}^{7}(q^{11};q^{11})_{\infty}^{3}}{q^5(q^2;q^2)_{\infty}^{3}(q^{22};q^{22})_{\infty}^{7}},\\ & \ \ \ \ \ \ \ \ \ \ \ \ \ \ \ \ \ \ \ \ \ \ \ \ \ \ \ \ \ \frac{5}{4}\frac{(q^2;q^2)_{\infty}(q^{11};q^{11})_{\infty}^{11}}{q^5(q;q)_{\infty}(q^{22};q^{22})_{\infty}^{11}}-\frac{3}{11}\frac{(q^2;q^2)_{\infty}^{8}(q^{11};q^{11})_{\infty}^{4}}{q^5(q;q)_{\infty}^{4}(q^{22};q^{22})_{\infty}^{8}} +\frac{1}{44}\frac{(q;q)_{\infty}^{7}(q^{11};q^{11})_{\infty}^{3}}{q^5(q^2;q^2)_{\infty}^{3}(q^{22};q^{22})_{\infty}^{7}}\}\\
&\{ \texttt{p}_{\texttt{g}}\texttt{(t):g}\in\texttt{AB}\}\texttt{: }  \ \ \ \{6776+9427 t+15477 t^2+13332 t^3+1078 t^4, -9581+594 t+5390 t^2+187 t^3,\\ & \ \ \ \ \ \ \ \ \ \ \ \ \ \ \ \ \ \ \ \ \ \ \ \ \ \ \ \ \ -6754+5368 t+2761 t^2+11 t^3 \}\\
&\texttt{Common Factor: } \ \ \ \ \ \ 11
\end{flalign*}

Our procedure returns a variation on a result that Radu already computed \cite{Radu}.  It serves as a witness identity for the divisibility of $p(11n+6)$ by 11, though it is not very satisfying.  It has a form resembling the classic witness identities which Ramanujan discovered for his congruences of $p(5n+4), p(7n+5)$ by $5, 7$, respectively.  In particular, the coefficients of $t$ in the membership witness are all divisible by 11.  Therefore, the result is a witness identity, provided one accepts that the functions of the algebra basis have integer power series expansions.  This is true, but not obvious.

In particular, we find a prevalence of 11 throughout the denominators of each function in our basis.  This is of course the one factor we would not want to find in the denominators!  Peter Paule was the first to realize that the integral expansions of the basis functions need to be proved, and successfully did so in \cite[Discussion, pp. 541-542]{Paule2}.

\subsection{An Identity for Broken 2-Diamond Partitions}

Broken $k$-diamond partitions, denoted by $\Delta_k(n)$, were defined by Andrews and Paule in 2007 \cite{AndP}.  They conjectured that

\begin{theorem}
For all $n\in\mathbb{Z}_{\ge 0}$,
\begin{align*}
\Delta_2(25n+14)\equiv \Delta_2(25n+24) \equiv 0\pmod{5}.
\end{align*}
\end{theorem}  This was subsequently proved in 2008 by Chan \cite{Chan}.  In 2015 Radu was able \cite{Radu} to give a proof by studying another arithmetic function with a simpler generating function.  Our complete implementation allows us to verify these congruences by directly examining the generating function for $\Delta_2(n)$.

We take $N=10, M=10, r=(-3,1,1,-1), m=25, j=14$.  Our package returns

\begin{flalign*}
\texttt{In[10] = }&\texttt{RK} [10,10,\{-3,1,1,-1\},25,14]&\\
\texttt{Out[10] = }&\\
&\texttt{P}_{\texttt{m,r}}\texttt{(j): } \ \ \ \ \ \ \ \ \ \ \ \ \ \ \ \  \{14, 24\}\\
&\texttt{f}_{\texttt{1}}\texttt{(q): } \ \ \ \ \ \ \ \ \ \ \ \ \ \ \ \ \ \ \  \frac{(q;q)_{\infty}^{126}(q^5;q^5)_{\infty}^{70}}{q^{58}(q^{2};q^{2})_{\infty}^{2}(q^{10};q^{10})_{\infty}^{190}}\\
&\texttt{t: }  \ \ \ \ \ \ \ \ \ \ \ \ \ \ \ \ \ \ \ \ \ \ \ \ \ \frac{(q^2;q^2)_{\infty}(q^5;q^5)_{\infty}^{5}}{q(q;q)_{\infty}(q^{10};q^{10})_{\infty}^{5}}\\
&\texttt{AB: } \ \ \ \ \ \ \ \ \ \ \ \ \ \ \ \ \ \ \ \ \ \ \ \ \{1\}\\
&\{ \texttt{p}_{\texttt{g}}\texttt{(t):g}\in\texttt{AB}\}\texttt{: }  \ \ \ \ \ \{...\}\\
&\texttt{Common Factor: } \ \ \ \ \ \ \ 25
\end{flalign*}

The membership witness returns a lengthy result, with terms of the order of $10^{76}$.  However, the computation time is short---less than 40 seconds with a 2.6 GHz Intel Processor on a modest laptop.  The complete witness is available, and easily computed, at \url{https://www3.risc.jku.at/people/nsmoot/RKAlg/RKSupplement1.nb}.

Each term in the membership witness is divisible by 25.  By expanding the generating function for $\Delta_2(n)$, one determines that $\Delta_2(14) = 10445$, and that $\Delta_2(49) = 1022063815$.

Because each of these numbers is divisible by 5 but not by 25, therefore $\sum_{n\ge 0} \Delta_2(25n+14)q^n$, $\sum_{n\ge 0} \Delta_2(25n+24)q^n$ must each be divisible by exactly one power of 5.  This completes the proof.

\subsection{Congruences with Overpartitions}

An enormous amount of work has been published in recent years on the congruence properties of overpartition functions, and our package has a great deal of utility in this subject.  We will examine three distinct problems here: two will involve the standard overpartition function $\bar{p}(n)$, and one will involve an overpartition function with additional restrictions $A_m(n)$.  In each case, we are able to make substantial improvements to previously established results.

As a preliminary, an overpartition of $n$ is a partition of $n$ in which the first occurrence of a part may or may not be ``marked."  Generally, this ``mark" is denoted with an overline (hence the term ``overpartition").  For example, the number 3 has 8 overpartitions:

\begin{align*}
&3,\\
&\bar{3},\\
&2+1,\\
&\bar{2}+1,\\
&2+\bar{1},\\
&\bar{2}+\bar{1},\\
&1+1+1,\\
&\bar{1}+1+1.
\end{align*}  We denote the number of overpartitions of $n$ by $\bar{p}(n)$.  The generating function for $\bar{p}(n)$ has the form

\begin{align*}
\sum_{n=0}^{\infty}\bar{p}(n)q^n = \frac{(-q;q)_{\infty}}{(q;q)_{\infty}} = \frac{(q^2;q^2)_{\infty}}{(q;q)^2_{\infty}}
\end{align*}  Part of the appeal of $\bar{p}(n)$ is the simplicity of the combinatoric interpretation, given the relative complexity of its generating function \cite{Corteel}.

\subsubsection{Congruences For $\bar{p}(n)$}

We will begin by giving some remarkable improvements to previously established congruences over $\bar{p}(n)$.  Moreover, we have the opportunity to apply our ``manual" procedure, and use the connection of modular functions with the topology of associated Riemann surfaces in order to construct a suitable algebra basis.

In 2016 Dou and Lin showed \cite{Dou} that

\begin{align}
\bar{p}(80n+8)\equiv \bar{p}(80n+52)\equiv \bar{p}(80n+68)\equiv \bar{p}(80n+72)\equiv 0\pmod{5}.\label{ovp8}
\end{align}  Hirschhorn in 2016 \cite{Hirschhorn1}, and Chern and Dastidar in 2018 \cite{Chern1} have studied these congruences as well, with the latter improving these congruences:

\begin{align*}
\bar{p}(80n+8)\equiv \bar{p}(80n+52)\equiv \bar{p}(80n+68)\equiv \bar{p}(80n+72)\equiv 0\pmod{25}.
\end{align*}  Chern and Dastidar go on to point out that

\begin{align*}
\bar{p}(135n+63)\equiv \bar{p}(135n+117)\equiv 0\pmod{5}.
\end{align*}

However, a quick computation of each of these sequences of overpartition numbers reveals much more.  For instance,

\begin{align*}
n\ \ \ & \bar{p}(80n+8)\\
0\ \ \ & 100\\
1\ \ \ & 8638130600\\
2\ \ \ & 350865646632400\\
3\ \ \ & 1512900775311002400\\
4\ \ \ & 1919738036947929590800\\
5\ \ \ & 1092453314947897908542800\\
6\ \ \ & 348534368588210202093102600\\
7\ \ \ & 71377855377904690816918291600\\
8\ \ \ & 10261762697785410674339371853700\\
\end{align*}

A very much stronger congruence clearly suggests itself.  We are able to make the following substantial improvements in each case:

\begin{theorem}
\begin{align*}
\bar{p}(80n+8)\equiv \bar{p}(80n+72)\equiv 0\pmod{100},\\
\bar{p}(80n+52)\equiv \bar{p}(80n+68)\equiv 0\pmod{200}.
\end{align*}
\end{theorem}

\begin{theorem}
\begin{align*}
\bar{p}(135n+63)\equiv \bar{p}(135n+117)\equiv 0\pmod{40}.
\end{align*}
\end{theorem}

Our package can be used to demonstrate each of these, though with some adjustments.  In the case of $\bar{p}(80n+j)$, we are forced to work over the congruence subgroup $\Gamma_0(40)$.  The generating set $\mathcal{G}_0(40)$ of the corresponding monoid $\mathcal{E}^{\infty}(40)$ of monopolar eta quotients can be computed with relative ease using \texttt{etaGenerators}; however, the set is nevertheless extremely large, and our procedure to compute the algebra basis using \texttt{AB} would be extremely inefficient.

We can remedy the problem by taking advantage of the Weierstrass gap theorem, (see \cite[Part 2, Section 17]{Weyl} for a classical introduction to the subject; see \cite{Paule3} for a more modern treatment of the theorem).  We use \cite[Theorem 3.1.1]{Diamond} to compute the genus of the corresponding modular curve $X_0(40)$ as 3, which implies that all modular functions with a pole only at $\infty$ on $\Gamma_0(40)$ must have order 4 or greater.  Radu's refinement of Newmann's conjecture \cite[Conjecture 9.4]{Paule1} suggests that a suitable combination of eta quotients will yield functions in $\left<\mathcal{E}^{\infty}(40)\right>_{\mathbb{Q}}$ with orders 4, 5, 6, 7.  Such a set of functions would be a sufficient algebra basis for $\left<\mathcal{E}^{\infty}(40)\right>_{\mathbb{Q}}$.

In this case, we are lucky, because a simple ordering of $\mathcal{G}_0(40)$ by the order of the elements at $\infty$ reveals that

\begin{align*}
\mathcal{G}_0(40)[1] = &\frac{(q^{4};q^{4})_{\infty}^{3}(q^{20};q^{20})_{\infty}}{q^4(q^{8};q^{8})_{\infty}(q^{40};q^{40})_{\infty}^{3}},\\
\mathcal{G}_0(40)[4] = &\frac{(q^{2};q^{2})_{\infty}^3(q^{5};q^{5})_{\infty}(q^{20};q^{20})_{\infty}^2}{q^5(q;q)_{\infty}(q^{10};q^{10})_{\infty}(q^{40};q^{40})_{\infty}^4},\\
\mathcal{G}_0(40)[7] = &\frac{(q^{2};q^{2})_{\infty}^6(q^{5};q^{5})_{\infty}^2(q^{8};q^{8})_{\infty}(q^{20};q^{20})_{\infty}^3}{q^6(q;q)_{\infty}^2(q^{4};q^{4})_{\infty}^3(q^{10};q^{10})_{\infty}^2(q^{40};q^{40})_{\infty}^5},\\
\mathcal{G}_0(40)[17] = &\frac{(q;q)_{\infty}^2(q^{5};q^{5})_{\infty}^2(q^{8};q^{8})_{\infty}^2(q^{20};q^{20})_{\infty}^3}{q^7(q^2;q^2)_{\infty}(q^{4};q^{4})_{\infty}(q^{10};q^{10})_{\infty}(q^{40};q^{40})_{\infty}^6}.
\end{align*}  Here for any ordered, enumerable set $\mathcal{S}$, we define the term $\mathcal{S}[j]$ as the $j$th term in the ordering of $\mathcal{S}$.

We can then define our algebra basis as 

\begin{align*}
T &= \mathcal{G}_0(40)[1],\\
\mathrm{Ab40} &= \left\{T,\left\{1, \mathcal{G}_0(40)[4], \mathcal{G}_0(40)[7], \mathcal{G}_0(40)[17]\right\} \right\}.
\end{align*}

Since we computed our algebra basis separately, we may now employ the manual case of our package, \texttt{RKMan} (See Section 2.4.2):

\begin{flalign*}
\texttt{In[11] = }&\texttt{RKMan} [40,2,\{-2,1\},80,8,\texttt{Ab40}]&\\
\texttt{Out[11] = }&\\
&\texttt{P}_{\texttt{m,r}}\texttt{(j): } \ \ \ \ \ \ \ \ \ \ \ \ \ \ \ \  \{8,72\}\\
&\texttt{f}_{\texttt{1}}\texttt{(q): } \ \ \ \ \ \ \ \ \ \ \ \ \ \ \ \ \ \ \  \frac{(q;q)_{\infty}^{333}(q^{8};q^{8})_{\infty}^{66}(q^{10};q^{10})_{\infty}^{36}(q^{20};q^{20})_{\infty}^{165}}{q^{400}(q^{2};q^{2})_{\infty}^{168}(q^{4};q^{4})_{\infty}^{31}(q^{5};q^{5})_{\infty}^{65}(q^{40};q^{40})_{\infty}^{334}}\\
&\texttt{t: }  \ \ \ \ \ \ \ \ \ \ \ \ \ \ \ \ \ \ \ \ \ \ \ \ \ \frac{(q^{4};q^{4})_{\infty}^{3}(q^{20};q^{20})_{\infty}}{q^4(q^{8};q^{8})_{\infty}(q^{40};q^{40})_{\infty}^{3}}\\
&\texttt{AB: } \ \ \ \ \ \ \ \ \ \ \ \ \ \ \ \ \ \ \ \ \ \ \ \ \{1, \frac{(q^{2};q^{2})_{\infty}^3(q^{5};q^{5})_{\infty}(q^{20};q^{20})_{\infty}^2}{q^5(q;q)_{\infty}(q^{10};q^{10})_{\infty}(q^{40};q^{40})_{\infty}^4}, \frac{(q^{2};q^{2})_{\infty}^6(q^{5};q^{5})_{\infty}^2(q^{8};q^{8})_{\infty}(q^{20};q^{20})_{\infty}^3}{q^6(q;q)_{\infty}^2(q^{4};q^{4})_{\infty}^3(q^{10};q^{10})_{\infty}^2(q^{40};q^{40})_{\infty}^5},\\ & \ \ \ \ \ \ \ \ \ \ \ \ \ \ \ \ \ \ \ \ \ \ \ \ \ \ \ \ \ \ \ \ \ \frac{(q;q)_{\infty}^2(q^{5};q^{5})_{\infty}^2(q^{8};q^{8})_{\infty}^2(q^{20};q^{20})_{\infty}^3}{q^7(q^2;q^2)_{\infty}(q^{4};q^{4})_{\infty}(q^{10};q^{10})_{\infty}(q^{40};q^{40})_{\infty}^6}\}\\
&\{ \texttt{p}_{\texttt{g}}\texttt{(t):g}\in\texttt{AB}\}\texttt{: }  \ \ \ \ \ \{...\}\\
&\texttt{Common Factor: } \ \ \ \ \ \ \ 10000
\end{flalign*}

The membership witness is too lengthy to present in this article.  The complete output of the algorithm can be found at \url{https://www3.risc.jku.at/people/nsmoot/RKAlg/RKSupplement2.nb}.  It is trivial to compute $\bar{p}(80n+~8),$ $\bar{p}(80n+72)$ for a handful of small $n$ in order to demonstrate that neither is divisible by $2^3$ or $5^3$.  Since the left hand side consists of a prefactor (with initial coefficient 1) and a product of the form

\begin{align*}
\left(\sum_{n\ge 0}\bar{p}(80n+8)q^n\right)\left(\sum_{n\ge 0}\bar{p}(80n+72)q^n \right),
\end{align*} with neither factor divisible by $2^3$ or $5^3$, the only remaining possibility is that each factor is divisible by $2^2\cdot 5^2 = 100$.

An almost identical output is produced for

\begin{flalign*}
\texttt{In[11] = }&\texttt{RKMan} [40,2,\{-2,1\},80,52,\texttt{Ab40}]&
\end{flalign*} but with an output of 40000 for congruences.  This is also available at \url{https://www3.risc.jku.at/people/nsmoot/RKAlg/RKSupplement2.nb}.  We may show that $\bar{p}(80n+52),$ $\bar{p}(80n+~68)$ are each divisible by 200, in a similar manner to the case of $\bar{p}(80n+8),$ $\bar{p}(80n+~72)$.

Finally, we consider the case of $\bar{p}(135n+63),\bar{p}(135n+117).$  We may similarly construct an algebra basis manually.  In this case, the most convenient congruence subgroup to work over is $\Gamma_0(30)$ ($N=30$).  The genus of $X_0(30)$ is 3, but we are at a slight disadvantage: there are eta quotients in $\mathcal{E}^{\infty}(30)$ with orders 4, 6, and 7, but none with order 5.  But we can construct a difference of eta quotients, each with order 6, to produce a function of order 5.  If we order the generators of $\mathcal{E}^{\infty}(30)$ by order at $\infty$, then

\begin{align*}
\mathcal{G}_0(30)[1] = &\frac{(q;q)_{\infty}(q^{6};q^{6})_{\infty}^{6}(q^{10};q^{10})_{\infty}^{2}(q^{15};q^{15})_{\infty}^{3}}{q^4(q^{2};q^{2})_{\infty}^2(q^{3};q^{3})_{\infty}^{3}(q^{5};q^{5})_{\infty}(q^{30};q^{30})_{\infty}^{6}},\\
\mathcal{G}_0(30)[4] - \mathcal{G}_0(30)[3] = &\frac{(q^2;q^2)_{\infty}^4(q^{10};q^{10})_{\infty}^{4}(q^{15};q^{15})_{\infty}^{4}}{q^6(q;q)_{\infty}^2(q^{5};q^{5})_{\infty}^{2}(q^{30};q^{30})_{\infty}^{8}}\\
&-\frac{(q;q)_{\infty}(q^{6};q^{6})_{\infty}^{2}(q^{10};q^{10})_{\infty}^{10}(q^{15};q^{15})_{\infty}^{5}}{q^6(q^{2};q^{2})_{\infty}^2(q^{3};q^{3})_{\infty}(q^{5};q^{5})_{\infty}^5(q^{30};q^{30})_{\infty}^{10}},
\end{align*}
\begin{align*}
\mathcal{G}_0(30)[2] = &\frac{(q;q)_{\infty}(q^{2};q^{2})_{\infty}(q^{5};q^{5})_{\infty}(q^{6};q^{6})_{\infty}(q^{10};q^{10})_{\infty}(q^{15};q^{15})_{\infty}^3}{q^6(q^{3};q^{3})_{\infty}(q^{30};q^{30})_{\infty}^{7}},\\
\mathcal{G}_0(30)[6] = &\frac{(q;q)_{\infty}(q^{5};q^{5})_{\infty}^{2}(q^{6};q^{6})_{\infty}(q^{10};q^{10})_{\infty}(q^{15};q^{15})_{\infty}^{3}}{q^7(q^{30};q^{30})_{\infty}^{8}}.
\end{align*} The orders here are (respectively) $4, 5, 6, 7$, again sufficient for an algebra basis:

\begin{align*}
T &= \mathcal{G}_0(30)[1],\\
G_1 &=\mathcal{G}_0(30)[4] - \mathcal{G}_0(30)[3]\\
G_2 &=\mathcal{G}_0(30)[2]\\
G_3 &=\mathcal{G}_0(30)[6]\\
\mathrm{Ab30} &= \left\{T,\left\{1, G_1, G_2, G_3\right\} \right\}.
\end{align*}  Employing \texttt{RKMan} once again, we get

\begin{flalign*}
\texttt{In[12] = }&\texttt{RKMan} [30,2,\{-2,1\},135,63,\texttt{Ab30}]&\\
\texttt{Out[12] = }&\\
&\texttt{P}_{\texttt{m,r}}\texttt{(j): } \ \ \ \ \ \ \ \ \ \ \ \ \ \ \ \  \{63, 117\}\\
&\texttt{f}_{\texttt{1}}\texttt{(q): } \ \ \ \ \ \ \ \ \ \ \ \ \ \ \ \ \ \ \  \frac{(q;q)_{\infty}^{653}(q^{6};q^{6})_{\infty}^{235}(q^{10};q^{10})_{\infty}^{272}(q^{15};q^{15})_{\infty}^{358}}{q^{507}(q^{2};q^{2})_{\infty}^{359}(q^{3};q^{3})_{\infty}^{275}(q^{5};q^{5})_{\infty}^{226}(q^{30};q^{30})_{\infty}^{656}}\\
&\texttt{t: }  \ \ \ \ \ \ \ \ \ \ \ \ \ \ \ \ \ \ \ \ \ \ \ \ \ \frac{(q;q)_{\infty}(q^{6};q^{6})_{\infty}^{6}(q^{10};q^{10})_{\infty}^{2}(q^{15};q^{15})_{\infty}^{3}}{q^4(q^{2};q^{2})_{\infty}^2(q^{3};q^{3})_{\infty}^{3}(q^{5};q^{5})_{\infty}(q^{30};q^{30})_{\infty}^{6}}\\
&\texttt{AB: } \ \ \ \ \ \ \ \ \ \ \ \ \ \ \ \ \ \ \ \ \ \ \ \ \{1, \frac{(q^2;q^2)_{\infty}^4(q^{10};q^{10})_{\infty}^{4}(q^{15};q^{15})_{\infty}^{4}}{q^6(q;q)_{\infty}^2(q^{5};q^{5})_{\infty}^{2}(q^{30};q^{30})_{\infty}^{8}} -\frac{(q;q)_{\infty}(q^{6};q^{6})_{\infty}^{2}(q^{10};q^{10})_{\infty}^{10}(q^{15};q^{15})_{\infty}^{5}}{q^6(q^{2};q^{2})_{\infty}^2(q^{3};q^{3})_{\infty}(q^{5};q^{5})_{\infty}^5(q^{30};q^{30})_{\infty}^{10}},\\
& \ \ \ \ \ \ \ \ \ \ \ \ \ \ \ \ \ \ \ \ \ \ \ \ \ \ \ \ \ \ \ \ \ \frac{(q;q)_{\infty}(q^{2};q^{2})_{\infty}(q^{5};q^{5})_{\infty}(q^{6};q^{6})_{\infty}(q^{10};q^{10})_{\infty}(q^{15};q^{15})_{\infty}^3}{q^6(q^{3};q^{3})_{\infty}(q^{30};q^{30})_{\infty}^{7}},\\
& \ \ \ \ \ \ \ \ \ \ \ \ \ \ \ \ \ \ \ \ \ \ \ \ \ \ \ \ \ \ \ \ \ \frac{(q;q)_{\infty}(q^{5};q^{5})_{\infty}^{2}(q^{6};q^{6})_{\infty}(q^{10};q^{10})_{\infty}(q^{15};q^{15})_{\infty}^{3}}{q^7(q^{30};q^{30})_{\infty}^{8}}\}\\
&\{ \texttt{p}_{\texttt{g}}\texttt{(t):g}\in\texttt{AB}\}\texttt{: }  \ \ \ \ \ \{...\}\\
&\texttt{Common Factor: } \ \ \ \ \ \ \ \frac{1600}{3}
\end{flalign*}

Once again, the membership witness is too large to present here.  It can be found in its entirety at \url{https://www3.risc.jku.at/people/nsmoot/RKAlg/RKSupplement2.nb}.  However, the fractional common factor emerges because each polynomial $p_{g}$ in the witness has integer coefficients, except for $p_{G_1}$, which is a polynomial in $\frac{1}{3}\mathbb{Z}$.  Because the remaining polynomials have integer coefficients (and all of the eta quotients involved have integer-coefficient expansions), we can conclude that $G_1$ has coefficients divisible by 3.  At any rate, this makes no difference for congruences with respect to powers of 2 or 5.

We may again quickly demonstrate that $\bar{p}(135n+63),\bar{p}(135n+117)$ are not divisible by $2^4$ or $5^2$, indicating that they must each be divisible by $2^3\cdot 5=40$.

\subsubsection{A Congruence for $\bar{p}(n)$ Modulo 243}

In 2017 Xia conjectured \cite{Xia} that

\begin{align*}
\bar{p}(96n+76)\equiv 0\pmod{3^5}
\end{align*} for all $n\in\mathbb{Z}_{\ge 0}$.  This conjecture was recently proved by Huang and Yao \cite{Huang}.

We have extended the theorem further:

\begin{theorem}
\begin{align*}
\bar{p}(96n+76)\equiv 0\pmod{2^3 3^5}
\end{align*} for all $n\in\mathbb{Z}_{\ge 0}$.
\end{theorem}

\begin{flalign*}
\texttt{In[13] = }&\texttt{RK} [24,2,\{-2,1\},96,76]&\\
\texttt{Out[13] = }&\\
&\texttt{P}_{\texttt{m,r}}\texttt{(j): } \ \ \ \ \ \ \ \ \ \ \ \ \ \ \ \  \{76\}\\
&\texttt{f}_{\texttt{1}}\texttt{(q): } \ \ \ \ \ \ \ \ \ \ \ \ \ \ \ \ \ \ \  \frac{(q;q)_{\infty}^{213}(q^{6};q^{6})_{\infty}^{33}(q^{8};q^{8})_{\infty}^{77}(q^{12};q^{12})_{\infty}^{113}}{q^{150}(q^{2};q^{2})_{\infty}^{107}(q^{3};q^{3})_{\infty}^{64}(q^{4};q^{4})_{\infty}^{37}(q^{24};q^{24})_{\infty}^{227}}\\
&\texttt{t: }  \ \ \ \ \ \ \ \ \ \ \ \ \ \ \ \ \ \ \ \ \ \ \ \ \ \frac{(q^{6};q^{6})_{\infty}^{3}(q^{8};q^{8})_{\infty}}{q^{2}(q^{2};q^{2})_{\infty}(q^{24};q^{24})_{\infty}^{3}}\\
&\texttt{AB: } \ \ \ \ \ \ \ \ \ \ \ \ \ \ \ \ \ \ \ \ \ \ \ \ \{1, \frac{(q^{6};q^{6})_{\infty}^{3}(q^{8};q^{8})_{\infty}}{q^{2}(q^{2};q^{2})_{\infty}(q^{24};q^{24})_{\infty}^{3}}+\frac{(q;q)_{\infty}(q^{3};q^{3})_{\infty}(q^{12};q^{12})_{\infty}(q^{4};q^{4})_{\infty}^{3}}{q^{3}(q^{2};q^{2})_{\infty}^{2}(q^{24};q^{24})_{\infty}^{4}}\}\\
&\{ \texttt{p}_{\texttt{g}}\texttt{(t):g}\in\texttt{AB}\}\texttt{: }  \ \ \ \ \ \{...\}\\
&\texttt{Common Factor: } \ \ \ \ \ \ \ 1944
\end{flalign*}

The theorem is then established, since $1944=2^3\cdot 3^5$.  The full identity can be found at \url{https://www3.risc.jku.at/people/nsmoot/RKAlg/RKSupplement2.nb}.

\subsubsection{A Restricted Overpartition Function}

Let $A_{m}(n)$ be the number of overpartitions of $n$ in which only the parts not divisible by $m$ may be overlined.  Then it can be showed that \cite{Munagi}

\begin{align*}
\sum_{n=0}^{\infty}A_m(n)q^n = \frac{(q^{2};q^{2})_{\infty}(q^m;q^m)_{\infty}}{(q;q)^2_{\infty}(q^{2m};q^{2m})_{\infty}}.
\end{align*}  In 2014, Munagi and Sellers give a variety of interesting congruences for $A_m(n)$.

For instance, \cite[Corollary 4.4, Theorem 4.5]{Munagi}:

\begin{theorem}
\begin{align*}
A_3(3n+1)&\equiv 0\pmod{2},\\
A_3(3n+2)&\equiv 0\pmod{4}.
\end{align*}
\end{theorem}  Both of these can be proved quickly with our package.  For example, to prove $A_3(3n+1)\equiv 0\pmod{2}$:

\begin{flalign*}
\texttt{In[14] = }&\texttt{RK} [6,6,\{-2,1,1,-1\},3,1]&\\
\texttt{Out[14] = }&\\
&\texttt{P}_{\texttt{m,r}}\texttt{(j): } \ \ \ \ \ \ \ \ \ \ \ \ \ \ \ \  \{1\}\\
&\texttt{f}_{\texttt{1}}\texttt{(q): } \ \ \ \ \ \ \ \ \ \ \ \ \ \ \ \ \ \ \  \frac{(q;q)_{\infty}^{3}(q^2;q^2)_{\infty}(q^3;q^3)_{\infty}^{6}}{q(q^6;q^6)_{\infty}^9}\\
&\texttt{t: }  \ \ \ \ \ \ \ \ \ \ \ \ \ \ \ \ \ \ \ \ \ \ \ \ \ \frac{(q;q)_{\infty}^5(q^3;q^3)_{\infty}}{q(q^2;q^2)_{\infty}(q^6;q^6)_{\infty}^5}\\
&\texttt{AB: } \ \ \ \ \ \ \ \ \ \ \ \ \ \ \ \ \ \ \ \ \ \ \ \ \{1\}\\
&\{ \texttt{p}_{\texttt{g}}\texttt{(t):g}\in\texttt{AB}\}\texttt{: }  \ \ \ \ \ \{16+2t\}\\
&\texttt{Common Factor: } \ \ \ \ \ \ \ 2
\end{flalign*}

On the other hand, \cite[Theorem 4.7, Theorem 4.9]{Munagi} $A_3(27n+26)\equiv 0\pmod{3}$, and $A_9(27n+24)\equiv 0\pmod{3}$.  Using our package, we can prove more:
\begin{theorem}
\begin{align*}
A_3(27n+26)&\equiv 0\pmod{12},\\
A_9(27n+24)&\equiv 0\pmod{24}.
\end{align*}
\end{theorem}

For example, to show that $A_9(27n+24)\equiv 0\pmod{24}$:

\begin{flalign*}
\texttt{In[15] = }&\texttt{RK} [6,18,\{-2,1,0,0,1,-1\},27,24]&\\
\texttt{Out[15] = }&\\
&\texttt{P}_{\texttt{m,r}}\texttt{(j): } \ \ \ \ \ \ \ \ \ \ \ \ \ \ \ \  \{24\}\\
&\texttt{f}_{\texttt{1}}\texttt{(q): } \ \ \ \ \ \ \ \ \ \ \ \ \ \ \ \ \ \ \  \frac{(q;q)_{\infty}^{47}(q^3;q^3)_{\infty}^{12}}{q^9(q^2;q^2)_{\infty}^7(q^6;q^6)_{\infty}^{51}}\\
&\texttt{t: }  \ \ \ \ \ \ \ \ \ \ \ \ \ \ \ \ \ \ \ \ \ \ \ \ \ \frac{(q;q)_{\infty}^5(q^3;q^3)_{\infty}}{q(q^2;q^2)_{\infty}(q^6;q^6)_{\infty}^5}\\
&\texttt{AB: } \ \ \ \ \ \ \ \ \ \ \ \ \ \ \ \ \ \ \ \ \ \ \ \ \{1\}\\
&\{ \texttt{p}_{\texttt{g}}\texttt{(t):g}\in\texttt{AB}\}\texttt{: }  \ \ \ \ \ \{7703510787293184+5456653474332672 t +1649478582927360 t^2\\& \ \ \ \ \ \ \ \ \ \ \ \ \ \ \ \ \ \ \ \ \ \ \ \ \ \ \ \ \ \ \ +276646783352832 t^3 +27989228519424 t^4+1735943602176 t^5\\& \ \ \ \ \ \ \ \ \ \ \ \ \ \ \ \ \ \ \ \ \ \ \ \ \ \ \ \ \ \ \ +63885293568 t^6+1269340416 t^7+10941888 t^8+22056 t^9\}\\
&\texttt{Common Factor: } \ \ \ \ \ \ \ 24
\end{flalign*}

We expect that a very large variety of other congruences and associated results for overpartition functions still await discovery, and that our package will prove extremely useful.

\subsection{Some Identities by Baruah and Sarmah}

For $r\in\mathbb{Z}$, define

\begin{align*}
\sum_{n = 0}^{\infty} p_r(n)q^n = (q;q)^r_{\infty}.
\end{align*}  In 2013 Baruah and Sarmah \cite{Baruah} gave a large variety of results for $p_r(n)$, all of which are accessible through our package.  One especially interesting example, \cite[Theorem 2.1, (2.10)]{Baruah} is not a congruence, but rather a simple identity:

\begin{theorem}
\begin{align*}
p_8(3n+1)&= 0.
\end{align*}
\end{theorem}

We can verify this by taking $M=1, r=(8), m=4, j=3, N=4$:

\begin{flalign*}
\texttt{In[16] = }&\texttt{RK} [4,1,\{8\},4,3]&\\
\texttt{Out[16] = }&\\
&\texttt{P}_{\texttt{m,r}}\texttt{(j): } \ \ \ \ \ \ \ \ \ \ \ \ \ \ \ \  \{3\}\\
&\texttt{f}_{\texttt{1}}\texttt{(q): } \ \ \ \ \ \ \ \ \ \ \ \ \ \ \ \ \ \ \  \frac{(q^2;q^2)_{\infty}^{12}}{q(q;q)_{\infty}^4(q^4;q^4)_{\infty}^{16}}\\
&\texttt{t: }  \ \ \ \ \ \ \ \ \ \ \ \ \ \ \ \ \ \ \ \ \ \ \ \ \ \frac{(q;q)_{\infty}^8}{q(q^4;q^4)_{\infty}^8}\\
&\texttt{AB: } \ \ \ \ \ \ \ \ \ \ \ \ \ \ \ \ \ \ \ \ \ \ \ \ \{1\}\\
&\{ \texttt{p}_{\texttt{g}}\texttt{(t):g}\in\texttt{AB}\}\texttt{: }  \ \ \ \ \ \{0\}\\
&\texttt{Common Factor: } \ \ \ \ \ \ \ 0
\end{flalign*}

Baruah and Sarmah list several congruences \cite[Theorem 5.1]{Baruah} which may easily be proved.  For example:

\begin{theorem}
\begin{align*}
p_{-4}(4n+3)\equiv &0\pmod{8},\\
p_{-8}(4n+3)\equiv &0\pmod{64},\\
p_{-2}(5n+2)\equiv p_{-2}(5n+3)\equiv p_{-2}(5n+4) \equiv &0\pmod{5},\\
p_{-4}(5n+3)\equiv p_{-4}(5n+4)\equiv &0\pmod{5}.
\end{align*}
\end{theorem}

We prove the first case by setting $M=1, r=(-4), m=4, j=3, N=8$.

\begin{flalign*}
\texttt{In[17] = }&\texttt{RK} [8,1,\{-4\},4,3]&\\
\texttt{Out[17] = }&\\
&\texttt{P}_{\texttt{m,r}}\texttt{(j): } \ \ \ \ \ \ \ \ \ \ \ \ \ \ \ \  \{3\}\\
&\texttt{f}_{\texttt{1}}\texttt{(q): } \ \ \ \ \ \ \ \ \ \ \ \ \ \ \ \ \ \ \  \frac{(q;q)_{\infty}^{19}(q^4;q^4)_{\infty}^{15}}{q^4(q^2;q^2)_{\infty}^8(q^8;q^8)_{\infty}^{22}}\\
&\texttt{t: }  \ \ \ \ \ \ \ \ \ \ \ \ \ \ \ \ \ \ \ \ \ \ \ \ \ \frac{(q^4;q^4)_{\infty}^{12}}{q(q^2;q^2)_{\infty}^4(q^8;q^8)_{\infty}^8}\\
&\texttt{AB: } \ \ \ \ \ \ \ \ \ \ \ \ \ \ \ \ \ \ \ \ \ \ \ \ \{1\}\\
&\{ \texttt{p}_{\texttt{g}}\texttt{(t):g}\in\texttt{AB}\}\texttt{: }  \ \ \ \ \ \{512 t+1408 t^2+480 t^3+40 t^4\}\\
&\texttt{Common Factor: } \ \ \ \ \ \ \ 8
\end{flalign*}

The other cases of this theorem can be proved similarly.

In another example, they prove \cite[Theorem 5.1, (5.3)]{Baruah} that $p_{-8}(8n+7)\equiv 0\pmod{2^{9}}$, but we prove even more:

\begin{theorem}
\begin{align*}
p_{-8}(8n+7)&\equiv 0\pmod{2^{11}}.
\end{align*}
\end{theorem}
We set $N=4$:

\begin{flalign*}
\texttt{In[18] = }&\texttt{RK} [4,1,\{-8\},8,7]&\\
\texttt{Out[18] = }&\\
&\texttt{P}_{\texttt{m,r}}\texttt{(j): } \ \ \ \ \ \ \ \ \ \ \ \ \ \ \ \  \{7\}\\
&\texttt{f}_{\texttt{1}}\texttt{(q): } \ \ \ \ \ \ \ \ \ \ \ \ \ \ \ \ \ \ \  \frac{(q;q)_{\infty}^{84}}{q^8(q^2;q^2)_{\infty}^4(q^4;q^4)_{\infty}^{72}}\\
&\texttt{t: }  \ \ \ \ \ \ \ \ \ \ \ \ \ \ \ \ \ \ \ \ \ \ \ \ \ \frac{(q;q)_{\infty}^{8}}{q(q^4;q^4)_{\infty}^8}\\
&\texttt{AB: } \ \ \ \ \ \ \ \ \ \ \ \ \ \ \ \ \ \ \ \ \ \ \ \ \{1\}\\
&\{ \texttt{p}_{\texttt{g}}\texttt{(t):g}\in\texttt{AB}\}\texttt{: }  \ \ \ \ \ \{576460752303423488+162129586585337856 t +18718085951258624t^2\\& \ \ \ \ \ \ \ \ \ \ \ \ \ \ \ \ \ \ \ \ \ \ \ \ \ \ \ \ \ \ +1139094046375936 t^3 +38970385760256 t^4+737593524224 t^5\\& \ \ \ \ \ \ \ \ \ \ \ \ \ \ \ \ \ \ \ \ \ \ \ \ \ \ \ \ \ \ +7041187840 t^6 +27033600 t^7+22528 t^8\}\\
&\texttt{Common Factor: } \ \ \ \ \ \ \ 2048
\end{flalign*}

\subsection{5-Regular Bipartitions}

In 2016 Liuquan Wang developed \cite{Wang} a large class of interesting congruences for the 5-regular bipartition function $B_5(n)$, with the generator

\begin{align*}
\sum_{n=0}^{\infty}B_5(n)q^n = \frac{(q^5;q^5)^2_{\infty}}{(q;q)^2_{\infty}}.
\end{align*}  Among many results were the following:

\begin{align*}
B_5(4n+3)&\equiv 0\pmod{5},\\
B_5(5n+2)&\equiv B_5(5n+3)\equiv B_5(5n+4)\equiv 0\pmod{5},\\
B_5(20n+7)&\equiv B_5(20n+19)\equiv 0\pmod{25}.
\end{align*}  We are able to make the following improvements:

\begin{theorem}
\begin{align*}
B_5(4n+3)&\equiv 0\pmod{10},\\
B_5(5n+2)&\equiv B_5(5n+3)\equiv B_5(5n+4)\equiv 0\pmod{5},\\
B_5(20n+7)&\equiv B_5(20n+19)\equiv 0\pmod{100}.
\end{align*}
\end{theorem}

\begin{flalign*}
\texttt{In[19] = }&\texttt{RK} [20,5,\{-2,2\},4,3]&\\
\texttt{Out[19] = }&\\
&\texttt{P}_{\texttt{m,r}}\texttt{(j): } \ \ \ \ \ \ \ \ \ \ \ \ \ \ \ \  \{3\}\\
&\texttt{f}_{\texttt{1}}\texttt{(q): } \ \ \ \ \ \ \ \ \ \ \ \ \ \ \ \ \ \ \ \frac{(q;q)_{\infty}^{6}(q^{2};q^{2})_{\infty}(q^{4};q^{4})_{\infty}(q^{10};q^{10})_{\infty}^{7}}{q^7(q^{5};q^{5})_{\infty}^{2}(q^{20};q^{20})_{\infty}^{13}}\\
&\texttt{t: }  \ \ \ \ \ \ \ \ \ \ \ \ \ \ \ \ \ \ \ \ \ \ \ \ \ \frac{(q^{4};q^{4})_{\infty}^{4}(q^{10};q^{10})_{\infty}^{2}}{q^2(q^{2};q^{2})_{\infty}^{2}(q^{20};q^{20})_{\infty}^{4}}\\
&\texttt{AB: } \ \ \ \ \ \ \ \ \ \ \ \ \ \ \ \ \ \ \ \ \ \ \ \ \{1, \frac{(q^{4};q^{4})_{\infty}(q^{5};q^{5})_{\infty}^{5}}{q^3(q;q)_{\infty}(q^{20};q^{20})_{\infty}^{5}} - \frac{(q^{4};q^{4})_{\infty}^{4}(q^{10};q^{10})_{\infty}^{2}}{q^2(q^{2};q^{2})_{\infty}^{2}(q^{20};q^{20})_{\infty}^{4}}\}\\
&\{ \texttt{p}_{\texttt{g}}\texttt{(t):g}\in\texttt{AB}\}\texttt{: }  \ \ \ \ \ \{50-40 t-50 t^2+40 t^3,-50+40 t+10 t^2\}\\
&\texttt{Common Factor: } \ \ \ \ \ \ \ 10
\end{flalign*}

\begin{flalign*}
\texttt{In[20] = }&\texttt{RK} [5,5,\{-2,2\},5,2]&\\
\texttt{Out[20] = }&\\
&\texttt{P}_{\texttt{m,r}}\texttt{(j): } \ \ \ \ \ \ \ \ \ \ \ \ \ \ \ \  \{2,4\}\\
&\texttt{f}_{\texttt{1}}\texttt{(q): } \ \ \ \ \ \ \ \ \ \ \ \ \ \ \ \ \ \ \ \frac{(q;q)_{\infty}^{20}}{q^2(q^{5};q^{5})_{\infty}^{20}}\\
&\texttt{t: }  \ \ \ \ \ \ \ \ \ \ \ \ \ \ \ \ \ \ \ \ \ \ \ \ \ \frac{((q;q)_{\infty})^6}{q((q^5;q^5)_{\infty})^6}\\
&\texttt{AB: } \ \ \ \ \ \ \ \ \ \ \ \ \ \ \ \ \ \ \ \ \ \ \ \ \{1\}\\
&\{ \texttt{p}_{\texttt{g}}\texttt{(t):g}\in\texttt{AB}\}\texttt{: }  \ \ \ \ \ \{15625+2500 t+100 t^2\}\\
&\texttt{Common Factor: } \ \ \ \ \ \ \ 25
\end{flalign*}

\begin{flalign*}
\texttt{In[21] = }&\texttt{RK} [5,5,\{-2,2\},5,3]&\\
\texttt{Out[21] = }&\\
&\texttt{P}_{\texttt{m,r}}\texttt{(j): } \ \ \ \ \ \ \ \ \ \ \ \ \ \ \ \  \{3\}\\
&\texttt{f}_{\texttt{1}}\texttt{(q): } \ \ \ \ \ \ \ \ \ \ \ \ \ \ \ \ \ \ \ \frac{(q;q)_{\infty}^{10}}{q(q^{5};q^{5})_{\infty}^{10}}\\
&\texttt{t: }  \ \ \ \ \ \ \ \ \ \ \ \ \ \ \ \ \ \ \ \ \ \ \ \ \ \frac{((q;q)_{\infty})^6}{q((q^5;q^5)_{\infty})^6}\\
&\texttt{AB: } \ \ \ \ \ \ \ \ \ \ \ \ \ \ \ \ \ \ \ \ \ \ \ \ \{1\}\\
&\{ \texttt{p}_{\texttt{g}}\texttt{(t):g}\in\texttt{AB}\}\texttt{: }  \ \ \ \ \ \{125+10 t\}\\
&\texttt{Common Factor: } \ \ \ \ \ \ \ 5
\end{flalign*}

\begin{flalign*}
\texttt{In[22] = }&\texttt{RK} [10,5,\{-2,2\},20,7]&\\
\texttt{Out[22] = }&\\
&\texttt{P}_{\texttt{m,r}}\texttt{(j): } \ \ \ \ \ \ \ \ \ \ \ \ \ \ \ \  \{7,19\}\\
&\texttt{f}_{\texttt{1}}\texttt{(q): } \ \ \ \ \ \ \ \ \ \ \ \ \ \ \ \ \ \ \ \frac{(q;q)_{\infty}^{77}(q^5;q^5)_{\infty}^{31}}{q^{27}(q^{2};q^{2})_{\infty}^{21}(q^{10};q^{10})_{\infty}^{87}}\\
&\texttt{t: }  \ \ \ \ \ \ \ \ \ \ \ \ \ \ \ \ \ \ \ \ \ \ \ \ \ \frac{(q^2;q^2)_{\infty}(q^5;q^5)_{\infty}^5}{q(q;q)_{\infty}(q^{10};q^{10})_{\infty}^5}\\
&\texttt{AB: } \ \ \ \ \ \ \ \ \ \ \ \ \ \ \ \ \ \ \ \ \ \ \ \ \{1\}\\
&\{ \texttt{p}_{\texttt{g}}\texttt{(t):g}\in\texttt{AB}\}\texttt{: }  \ \ \ \ \ \{7388718138654720000 t^2+153008038121308160000 t^3\\& \ \ \ \ \ \ \ \ \ \ \ \ \ \ \ \ \ \ \ \ \ \ \ \ +1257731351012966400000 t^4+5675499664745431040000 t^5\\& \ \ \ \ \ \ \ \ \ \ \ \ \ \ \ \ \ \ \ \ \ \ \ \ +16507857641427435520000 t^6+34080767872618987520000 t^7\\& \ \ \ \ \ \ \ \ \ \ \ \ \ \ \ \ \ \ \ \ \ \ \ \ +53266856094927421440000 t^8+65937188949118156800000 t^9\\& \ \ \ \ \ \ \ \ \ \ \ \ \ \ \ \ \ \ \ \ \ \ \ \ +66700597538020392960000 t^{10}+56314162511641313280000 t^{11}\\& \ \ \ \ \ \ \ \ \ \ \ \ \ \ \ \ \ \ \ \ \ \ \ \ +40234227634725191680000 t^{12}+24527816166851215360000 t^{13}
\\& \ \ \ \ \ \ \ \ \ \ \ \ \ \ \ \ \ \ \ \ \ \ \ \ +12802067441385472000000 t^{14}+5714660420762992640000 t^{15}
\\& \ \ \ \ \ \ \ \ \ \ \ \ \ \ \ \ \ \ \ \ \ \ \ \ +2169098785981726720000 t^{16}+691839480120197120000 t^{17}
\\& \ \ \ \ \ \ \ \ \ \ \ \ \ \ \ \ \ \ \ \ \ \ \ \ +181850756413399040000 t^{18}+38175700204339200000 t^{19}
\\& \ \ \ \ \ \ \ \ \ \ \ \ \ \ \ \ \ \ \ \ \ \ \ \ +6075890734530560000 t^{20}+680092466755680000 t^{21}
\\& \ \ \ \ \ \ \ \ \ \ \ \ \ \ \ \ \ \ \ \ \ \ \ \ +49080942745680000 t^{22}+2083485921960000 t^{23}+46908276350000 t^{24}
\\& \ \ \ \ \ \ \ \ \ \ \ \ \ \ \ \ \ \ \ \ \ \ \ \ +483406090000 t^{25}+1812970000 t^{26}+1190000 t^{27}\}\\
&\texttt{Common Factor: } \ \ \ \ \ \ \ 10000
\end{flalign*}

\subsection{Some Congruences Related to the Tau Function}

Please note that in this section we will assume $q=e^{2\pi i z}$ with $z\in\mathbb{H}$, to avoid confusion with $\tau$, which will be used to identify a certain arithmetic function.

Ramanujan's tau function is defined by the following:

\begin{align*}
\Delta(z) := \sum_{n=1}^{\infty}\tau(n) q^n = q(q;q)^{24} = \eta(z)^{24}.
\end{align*}  The functions $\Delta(z)$ and $\tau(n)$ are among the most studied objects in the theory of modular forms.  In particular, numerous interesting congruences have been found.  Many classic examples include the following, discovered by Ramanujan \cite{Ramanujan2}:

\begin{theorem}
\begin{align*}
\tau(7n+m) \equiv 0\pmod{7}
\end{align*} for $m\in\{0, 3, 5, 6\}$.
\end{theorem}  Our algorithm can easily handle each of these cases.  For example, we take the case of $\tau(7n)$ (notice that we study $(q;q)_{\infty}^{24}$, rather than $q(q;q)_{\infty}^{24}$; because of this, we need to examine the progression $7n+6$):

\begin{flalign*}
\texttt{In[23] = }&\texttt{RK} [7,1,\{24\},7,6]&\\
\texttt{Out[23] = }&\\
&\texttt{P}_{\texttt{m,r}}\texttt{(j): } \ \ \ \ \ \ \ \ \ \ \ \ \ \ \ \  \{6\}\\
&\texttt{f}_{\texttt{1}}\texttt{(q): } \ \ \ \ \ \ \ \ \ \ \ \ \ \ \ \ \ \ \ \frac{1}{q^6(q^7;q^7)_{\infty}^{24}}\\
&\texttt{t: }  \ \ \ \ \ \ \ \ \ \ \ \ \ \ \ \ \ \ \ \ \ \ \ \ \ \frac{(q;q)_{\infty}^4}{q(q^7;q^7)_{\infty}^4}\\
&\texttt{AB: } \ \ \ \ \ \ \ \ \ \ \ \ \ \ \ \ \ \ \ \ \ \ \ \ \{1\}\\
&\{ \texttt{p}_{\texttt{g}}\texttt{(t):g}\in\texttt{AB}\}\texttt{: }  \ \ \ \ \ \{-1977326743-16744 t^6\}\\
&\texttt{Common Factor: } \ \ \ \ \ \ \ 7
\end{flalign*}

We will give a more recent example discovered by Koustav Banerjee \cite{Banerjee}:

\begin{theorem}
\begin{align*}
\tau(8(14n+k))\equiv 0\pmod{2^3\cdot 3\cdot 5\cdot 11},
\end{align*} for all $n\in\mathbb{Z}_{\ge 0}$ and $k$ an odd integer mod 14.
\end{theorem}  This may be broken up into three distinct RK identities.  We give the case of $112n+56$ (here shifted to $112n+55$)

\begin{flalign*}
\texttt{In[24] = }&\texttt{RK} [14,1,\{24\},112,55]&\\
\texttt{Out[24] = }&\\
&\texttt{P}_{\texttt{m,r}}\texttt{(j): } \ \ \ \ \ \ \ \ \ \ \ \ \ \ \ \  \{55\}\\
&\texttt{f}_{\texttt{1}}\texttt{(q): } \ \ \ \ \ \ \ \ \ \ \ \ \ \ \ \ \ \ \ \frac{(q^{2};q^{2})_{\infty}^{12}(q^{7};q^{7})_{\infty}^{30}}{q^{25}(q;q)_{\infty}^{6}(q^{14};q^{14})_{\infty}^{60}}\\
&\texttt{t: }  \ \ \ \ \ \ \ \ \ \ \ \ \ \ \ \ \ \ \ \ \ \ \ \ \ \frac{(q^{2};q^{2})_{\infty}(q^{7};q^{7})_{\infty}^{7}}{q^2(q;q)_{\infty}(q^{14};q^{14})_{\infty}^{7}}\\
&\texttt{AB: } \ \ \ \ \ \ \ \ \ \ \ \ \ \ \ \ \ \ \ \ \ \ \ \ \{1, \frac{(q^{2};q^{2})_{\infty}^8(q^{7};q^{7})_{\infty}^{4}}{q^3(q;q)_{\infty}^4(q^{14};q^{14})_{\infty}^{8}} - 4 \frac{(q^{2};q^{2})_{\infty}(q^{7};q^{7})_{\infty}^{7}}{q^2(q;q)_{\infty}(q^{14};q^{14})_{\infty}^{7}}\}\\
&\{ \texttt{p}_{\texttt{g}}\texttt{(t):g}\in\texttt{AB}\}\texttt{: }  \ \ \ \ \ \{1483245480837120+22804899267870720 t\\& \ \ \ \ \ \ \ \ \ \ \ \ \ \ \ \ \ \ \ \ \ \ \ \ \ \ \ \ \ \ \ -281353127146291200 t^2+4813307313059266560 t^3\\& \ \ \ \ \ \ \ \ \ \ \ \ \ \ \ \ \ \ \ \ \ \ \ \ \ \ \ \ \ \ \ -2117115491136307200 t^4-3347863578673152000 t^5 \\& \ \ \ \ \ \ \ \ \ \ \ \ \ \ \ \ \ \ \ \ \ \ \ \ \ \ \ \ \ \ \ +845098635118510080 t^6+77358598094131200 t^7\\& \ \ \ \ \ \ \ \ \ \ \ \ \ \ \ \ \ \ \ \ \ \ \ \ \ \ \ \ \ \ \ -25371836549283840 t^8-1132615297820160 t^9\\& \ \ \ \ \ \ \ \ \ \ \ \ \ \ \ \ \ \ \ \ \ \ \ \ \ \ \ \ \ \ \ -512964938787840 t^{10}-114993988032000 t^{11}-349389680640 t^{12},\\& \ \ \ \ \ \ \ \ \ \ \ \ \ \ \ \ \ \ \ \ \ \ \ \ \ \ \ \ \ \ \ -1483245480837120-6489198978662400 t+990900684041748480 t^2\\& \ \ \ \ \ \ \ \ \ \ \ \ \ \ \ \ \ \ \ \ \ \ \ \ \ \ \ \ \ \ \ -151791226737131520 t^3-1234180893392240640 t^4\\& \ \ \ \ \ \ \ \ \ \ \ \ \ \ \ \ \ \ \ \ \ \ \ \ \ \ \ \ \ \ \ +461934380423577600 t^5-65498418207129600 t^6\\& \ \ \ \ \ \ \ \ \ \ \ \ \ \ \ \ \ \ \ \ \ \ \ \ \ \ \ \ \ \ \ +2233732210913280 t^7+170807954042880 t^8+855016378191360 t^9\\& \ \ \ \ \ \ \ \ \ \ \ \ \ \ \ \ \ \ \ \ \ \ \ \ \ \ \ \ \ \ \ -4703322624000 t^{10}-1414533120 t^{11}\}\\
&\texttt{Common Factor: } \ \ \ \ \ \ \ 591360
\end{flalign*}

The congruence here is even stronger than in the more general case, since $591360 = 2^9\cdot 3\cdot 5\cdot 7\cdot 11$.

\subsection{An Identity Related to Rogers--Ramanujan Subpartitions}

We finish with an application of our package to studying infinite families of congruences.  In 2017 Choi, Kim, and Lovejoy discovered a congruence \cite[Proposition 6.4]{CKL}, based on a subpartition function studied by Kolitsch \cite{Kolitsch}.

For any partition $\lambda$, define the corresponding Rogers--Ramanujan subpartition of $\lambda$ as the unique subpartition of $\lambda$ with a maximal number of parts, in which the parts are nonrepeating, nonconsecutive, and larger than the remaining parts of $\lambda$.  For example, the partition $8+5+3+2+2+1+1+1$ contains the Rogers--Ramanujan subpartition $8+5+3$, whereas the partition $8+8+2+2+1+1+1$ contains the empty Rogers--Ramanujan subpartition.

Let us define $R_l(n)$ as the number of partitions of $n$ which contain a Rogers--Ramanujan subpartition of length $l$, and then

\begin{align*}
A(n) = \sum_{l\ge 0} l\cdot R_l(n).
\end{align*}

In \cite{CKL} the following was demonstrated:

\begin{theorem}
For $n\in\mathbb{Z}_{\ge 0}$,
\begin{align*}
A(25n+9)\equiv A(25n+14)\equiv A(25n+24)\equiv 0\pmod{5}.
\end{align*}
\end{theorem}  This was proved by connecting $A(n)$ with the coefficient $a(n)$ of

\begin{align*}
\sum_{n=0}^{\infty}a(n)q^n = \frac{(q^2;q^2)^5_{\infty}}{(q;q)^3_{\infty}(q^4;q^4)^2_{\infty}}.
\end{align*}   The authors of \cite{CKL} pointed to other suspected congruences and compared the generating function for $a(n)$ with that of $c\phi_2(n)$.  From this, they conjectured the existence of an infinite family for $A(n)$, in the style of Ramanujan's classic congruences, modulo powers of 5 \cite[Chapter 7]{Knopp}.  This infinite family was given a precise formulation after careful investigation using our standard package, as well as a modified version \cite{RadS} of the package designed to check large congruences.  After substantial evidence was gathered, the conjecture was proved \cite{Smoot}.  We will here consider the case $A(25n+24)$ by examining $a(25n+24)$.

Taking $M=4, r=(-3,5,-2), m=25, j=24$, and setting $N=20$, we find that

\begin{flalign*}
\texttt{In[25] = }&\texttt{RK} [20,4,\{-3,5,-2\},25,24]&\\
\texttt{Out[25] = }&\\
&\texttt{P}_{\texttt{m,r}}\texttt{(j): } \ \ \ \ \ \ \ \ \ \ \ \ \ \ \ \  \{24\}\\
&\texttt{f}_{\texttt{1}}\texttt{(q): } \ \ \ \ \ \ \ \ \ \ \ \ \ \ \ \ \ \ \ \frac{(q;q)_{\infty}^{35}(q^{4};q^{4})_{\infty}^{18}(q^{10};q^{10})_{\infty}^{30}}{q^{26}(q^{2};q^{2})_{\infty}^{27}(q^{5};q^{5})_{\infty}^{8}(q^{20};q^{20})_{\infty}^{48}}\\
&\texttt{t: }  \ \ \ \ \ \ \ \ \ \ \ \ \ \ \ \ \ \ \ \ \ \ \ \ \ \frac{(q^{4};q^{4})_{\infty}^{4}(q^{10};q^{10})_{\infty}^{2}}{q^2(q^{2};q^{2})_{\infty}^{2}(q^{20};q^{20})_{\infty}^{4}}\\
&\texttt{AB: } \ \ \ \ \ \ \ \ \ \ \ \ \ \ \ \ \ \ \ \ \ \ \ \ \{1, \frac{(q^{4};q^{4})_{\infty}(q^{5};q^{5})_{\infty}^{5}}{q^3(q;q)_{\infty}(q^{20};q^{20})_{\infty}^{5}} - \frac{(q^{4};q^{4})_{\infty}^{4}(q^{10};q^{10})_{\infty}^{2}}{q^2(q^{2};q^{2})_{\infty}^{2}(q^{20};q^{20})_{\infty}^{4}}\}\\
&\{ \texttt{p}_{\texttt{g}}\texttt{(t):g}\in\texttt{AB}\}\texttt{: }  \ \ \ \ \ \{126953125+74218750 t-174609375 t^2+25390625 t^3\\ & \ \ \ \ \ \ \ \ \ \ \ \ \ \ \ \ \ \ \ \ \ \ \ \ \ \ \ \ \ \ \ -1237031250 t^4+1542084375 t^5+3798876250 t^6\\ & \ \ \ \ \ \ \ \ \ \ \ \ \ \ \ \ \ \ \ \ \ \ \ \ \ \ \ \ \ \ \ -7568402750 t^7+3755535625 t^8+210440100 t^9\\ & \ \ \ \ \ \ \ \ \ \ \ \ \ \ \ \ \ \ \ \ \ \ \ \ \ \ \ \ \ \ \ -754603995 t^{10}+190492925 t^{11}+10649860 t^{12}+5735 t^{13},\\& \ \ \ \ \ \ \ \ \ \ \ \ \ \ \ \ \ \ \ \ \ \ \ \ \ \ \ \ \ \ \ -78125000+62500000 t-46093750 t^2+128906250 t^3\\& \ \ \ \ \ \ \ \ \ \ \ \ \ \ \ \ \ \ \ \ \ \ \ \ \ \ \ \ \ \ \ +551875000 t^4-1636475000 t^5+430767500 t^6+1615951500 t^7\\& \ \ \ \ \ \ \ \ \ \ \ \ \ \ \ \ \ \ \ \ \ \ \ \ \ \ \ \ \ \ \ -1247744000 t^8+145803400 t^9+72090170 t^{10}+543930 t^{11}\}\\
&\texttt{Common Factor: } \ \ \ \ \ \ \ 5
\end{flalign*}

\subsection{Other Examples}

The examples above are only a small sample of results produced in the last few years which our algorithm can duplicate or even improve on.  We have similarly applied our software to a large variety of other examples which, for simple want of space, we are unable to explicitly give.  Examples include identities found in \cite{Baruah2}, \cite{hchan}, \cite{Chen}, \cite{Fang}, \cite{Hirschhorn2}, \cite{Jameson}, \cite{Keith}, \cite{Lin}, \cite{Naika}, \cite{Tang}.

In summary, we have largely reduced the problem of discovering and proving such identities to a relatively fast computational procedure which can easily duplicate a large number of results that are frequently published in the mathematical literature.  We hope that researchers will make good use of our work.

\section{Accessibility}

Our software package is freely available as \texttt{RaduRK.m} via \url{https://www3.risc.jku.at/people/nsmoot/RKAlg/RaduRK.m}.  The implementation uses Mathematica, and requires installation of a Diophantine software package called \texttt{4ti2} \cite{4ti2}.  In particular, we used the interface \texttt{math4ti2.m} developed by Ralf Hemmecke and Silviu Radu.  We will also make our software available on the Computer Algebra for Combinatorics section of the RISC webpage \url{https://risc.jku.at/research_topic/computer-algebra-for-combinatorics/}.

A demonstration of the software can be found at \url{https://www3.risc.jku.at/people/nsmoot/RKAlg/RKSupplement1.nb}, in which most of the examples of this paper are computed; and \url{https://www3.risc.jku.at/people/nsmoot/RKAlg/RKSupplement2.nb}, in which the lengthier computations of Sections 3.5.1-3.5.2 are given.

Because \texttt{4ti2} is a Linux program, some additional steps are necessary in order to properly install our software onto an Apple or Windows operating system.  We provide the necessary steps for installation onto Apple and Windows at \url{https://www3.risc.jku.at/people/nsmoot/RKAlg/4ti2installationinstructions.rtf}.  All difficulties in installation should be communicated immediately to the author's email, \url{nsmoot@risc.uni-linz.ac.at}.

\section{Acknowledgments}

This research was funded by the Austrian Science Fund (FWF): W1214-N15, project DK6, and by the strategic program ``Innovatives OÖ 2010 Plus" by the Upper Austrian Government.  I wish to thank the Austrian government for its generous support.

My dearest thanks to each of the three anonymous referees, who supplied extremely useful criticisms.  I believe that my attempts to respond to them have substantially improved this paper.

I am extremely grateful to Drs. Ralf Hemmecke and Silviu Radu, who assisted me with their extremely useful software package \texttt{math4ti2.m} \cite{4ti22}.  In addition, I am deeply indebted to Ralf Wahner, who was exceedingly kind in his assistance with my installation instructions.  My software is available for installation onto an Apple, Windows, or Linux operating system.  This availability is owed almost entirely to Wahner's help.  RISC is very lucky to have his expertise. 

Finally, I am enormously grateful to Professor Peter Paule for recommending this remarkable project to me in the first place.  I am thankful for his guidance on technical matters, as well as his patience and kindness.


\begin{thebibliography}{X}

\bibitem{4ti2} 4ti2 team.  4ti2---A Software Package for Algebraic, Geometric and Combinatorial Problems on Linear Spaces: \url{<http://www.4ti2.de>}.

\bibitem{Andrews} G.E. Andrews, \textit{The Theory of Partitions}, \textit{Encyclopedia of Mathematics and its Applications} 2, Addison-Wesley (1976). Reissued, Cambridge (1998).

\bibitem{AndF} G. E. Andrews, ``Generalized Frobenius Partitions," \textit{Memoirs of the American Mathematical Society} 49, pp. 301 (1984).

\bibitem{AndP} G.E. Andrews, P. Paule, ``MacMahon's Partition Analysis XI: Broken Diamonds and Modular Forms," \textit{Acta Arithmetica} 126, pp. 281-294 (2007).

\bibitem{Banerjee} K. Banerjee, Private correspondence (February, 2019).

\bibitem{Baruah2} N.D. Baruah, N.M. Begum, ``Generating Functions and Congruences for Some Partition Functions Related to Mock Theta Functions'' (Submitted), Available at \url{https://arxiv.org/abs/1908.07741} (2019).

\bibitem{Baruah} N. D. Baruah, B. K. Sarmah, ``Identities and Congruences for the General Partition and Ramanujan's Tau Functions," \textit{Indian Journal of Pure and Applied Mathematics} 44 (5), pp. 643-671 (2013).

\bibitem{hchan} H. Chan, ``Ramanujan's Cubic Continued Fraction and an Analog of his `Most Beautiful Identity','' \textit{International Journal of Number Theory} 6 (3), pp. 673-680 (2010).

\bibitem{Chan} S. H. Chan.  ``Some Congruences for Andrews--Paule's Broken 2-Diamond Partitions," \textit{Discrete Mathematics} 308 (23), pp. 5735-5741 (2008).

\bibitem{Chen} W.Y.C. Chen, B.L.S. Lin, ``Congruences for the Number of Cubit Partitions Derived from Modular Forms,'' \textit{https://arxiv.org/abs/0910.1263} (2009).

\bibitem{Chern1} S. Chern, M. G. Dastidar, ``Some Congruences Modulo 5 and 25 for Overpartitions," \textit{Ramanujan Journal} 47 (2), pp. 435-445 (2018).

\bibitem{CKL} Y. Choi, B. Kim, J. Lovejoy, ``Overpartitions into Distinct Parts Without Short Sequences," \textit{Journal of Number Theory} 175, pp. 117-133 (2017).

\bibitem{Corteel} S. Corteel, J. Lovejoy, ``Overpartitions," \textit{Transactions of the American Mathematical Society} 356, pp. 1623-1635 (2004).

\bibitem{Diamond} F. Diamond, J. Shurman, \textit{A First Course in Modular Forms}, 4th Printing., Springer Publishing (2016).

\bibitem{Dou} D.Q.J. Dou, B.L.S. Lin, ``New Ramanujan-Type Congruences Modulo 5 for Overpartitions," \textit{Ramanujan Journal} 44 (2), pp. 401-410 (2016).

\bibitem{Euler} L. Euler, \textit{Introductio in Analysin Infinitorum, Chapter 16}, Marcum-Michaelem Bousquet, Lausannae (1748).

\bibitem{Fang} H. Fang, F. Xue, O.X.M. Yao, ``New Congruences Modulo 5 and 9 for Partitions with Odd Parts Distinct,'' \textit{Quaestiones Mathematicae}, pp. 1-14 (2019).

\bibitem{Hemmecke} R. Hemmecke, ``Dancing Samba with Ramanujan Partition Congruences," \textit{Journal of Symbolic Computation} 84, pp. 14-24 (2018).

\bibitem{4ti22} R. Hemmecke, S. Radu, 4ti2 Mathematical Interface: \url{<https://www3.risc.jku.at/research/combinat/software/math4ti2/math4ti2.m>}. (2017).

\bibitem{Hirschhorn1} M. Hirschhorn, ``Some Congruences for Overpartitions," \textit{New Zealand Journal of Mathematics} 46, pp. 141-144 (2016).

\bibitem{Hirschhorn2} M. Hirschhorn, ``A Conjecture of B. Lin on Cubic Partition Pairs,'' \textit{Ramanujan Journal} 45, pp. 781-795 (2018).

\bibitem{Huang} X. Huang, O.X.M. Yao, ``Proof of a Conjecture on a Congruence Modulo 243 for Overpartitions,'' \textit{Periodica Mathematica Hungarica} 79, pp. 227-235 (2019).

\bibitem{Jameson} M. Jameson, ``Congruences for Broken $k$-Diamond Partitions,'' \textit{Annals of Combinatorics} 17, pp. 333-338 (2013).

\bibitem{Keith} W.J. Keith, ``Restricted k-Color Partitions II'' (Submitted), Available at \url{https://arxiv.org/abs/2001.08351} (2020).

\bibitem{Knopp} M. Knopp, \textit{Modular Functions in Analytic Number Theory}, 2nd Ed., AMS Chelsea Publishing (1993).

\bibitem{Kolberg0} O. Kolberg, ``An Elementary Discussion of Certain Modular Forms,'' \textit{Univ. Bergen Arb. naturv. r.} 19 (1959).

\bibitem{Kolberg} O. Kolberg, ``Some Identities Involving the Partition Function," \textit{Mathematica Scandinavica} 5, pp. 77-92 (1957).

\bibitem{Kolitsch} L.W. Kolitsch, ``Rogers--Ramanujan Subpartitions and Their Connections to Other Partitions," \textit{Ramanujan Journal} 16 (2), pp. 163-167 (2008).

\bibitem{Lehner} J. Lehner, \textit{Discontinuous Groups and Automorphic Functions}, Mathematical Surveys and Monographs Number 8, American Mathematical Society (1964).

\bibitem{Lin} B.L.S. Lin, ``Congruences Modulo 27 for Cubic Partition Pairs,'' \textit{Journal of Number Theory} 171, pp. 31-42 (2017).

\bibitem{Miranda} R. Miranda, \textit{Algebraic Curves and Riemann Surfaces,} \textit{Grad. Stud. Math. Vol. 5}, AMS (1995).

\bibitem{Munagi} A. O. Munagi, J. Sellers, ``Refining Overlined Parts in Overpartitions via Residue Classes: Bijections, Generating Functions, and Congruences," \textit{Utilitas Mathematica} 95, pp. 33-49 (2014).

\bibitem{Naika} M.S.M. Naika, S.S. Nayaka, ``Congruences for Partition Quadruples with $t$-Cores,'' \textit{Acta Mathematica Vietnamica} (2019).

\bibitem{Newman1} M. Newman, ``Construction and Application of a Class of Modular Functions," \textit{Proceedings of the London Mathematical Society} 3 (7), pp. 334-350 (1956).

\bibitem{Newman2} M. Newman, ``Construction and Application of a Class of Modular Functions II," \textit{Proceedings of the London Mathematical Society} 3 (9), pp. 373-387 (1959).

\bibitem{Paule}  P. Paule, S. Radu, ``The Andrews--Sellers Family of Partition Congruences," \textit{Advances in Mathematics} 230 (3), pp. 819-838 (2012).

\bibitem{Paule2}  P. Paule, S. Radu, ``Partition Analysis, Modular Functions, and Computer Algebra," \textit{Recent Trends in Combinatorics, IMA Volume in Mathematics and its Applications, Springer}, pp. 511-543 (2015).

\bibitem{Paule3} P. Paule, S. Radu, ``A Proof of the Weierstrass Gap Theorem not using the Riemann--Roch Formula," \textit{Annals of Combinatorics} 23, pp. 963-1007 (2019).

\bibitem{Paule1}  P. Paule, S. Radu, ``A Unified Algorithmic Framework for Ramanujan's Congruences Modulo Powers of 5, 7, and 11," (Submitted), Available at \url{http://www3.risc.jku.at/publications/download/risc_5760/PP_submission.pdf} (2018).

\bibitem{Rademacher} H. Rademacher, ``The Ramanujan Identities Under Modular Substitutions,'' \textit{Transactions of the American Mathematical Society} 51 (3), pp. 609-636 (1942).

\bibitem{Radu0} S. Radu, ``An Algorithmic Approach to Ramanujan's Congruences," \textit{Ramanujan Journal} 20, pp. 215-251 (2009).

\bibitem{Radu} S. Radu, ``An Algorithmic Approach to Ramanujan--Kolberg Identities," \textit{Journal of Symbolic Computation} 68, pp. 225-253 (2015).

\bibitem{RadS} S. Radu, N. Smoot, ``A Method of Verifying Partition Congruences by Symbolic Computation," (Submitted), Available at \url{https://www3.risc.jku.at/publications/download/risc_5897/paper2b.pdf} (2019).

\bibitem{Ramanujan} S. Ramanujan, ``Some Properties of $p(n)$, the Number of Partitions of $n$", \textit{Proceedings of the Cambridge Philosophical Society}, pp. 207-210 (1919).

\bibitem{Ramanujan2} S. Ramanujan, ``Congruence Properties of Partitions", \textit{Proceedings of the London Mathematical Society}, 2 (18) (1920).

\bibitem{Sellers} J. Sellers, ``Congruences Involving $F$-Partition Functions," \textit{International Journal of Mathematics and Mathematical Sciences} 17, pp. 187-188 (1994).

\bibitem{Smoot} N. Smoot, ``A Family of Congruences for Rogers--Ramanujan Subpartitions," \textit{Journal of Number Theory} 196, pp. 35-60 (2019).

\bibitem{Tang} D. Tang, ``Congruences Modulo Powers of 3 for 2-Color Partition Triples,'' \textit{Periodica Mathematica Hungarica} 78, pp. 254-266 (2019).

\bibitem{Wang} L. Wang, ``Congruences Modulo Powers of 5 for Two Restricted Bipartitions," \textit{Ramanujan Journal} 44, pp. 471-491 (2014).

\bibitem{Watson} G.N. Watson, ``Ramanujans Vermutung über Zerfallungsanzahlen," \textit{J. Reine Angew. Math.,} pp. 97-118 (1938).

\bibitem{Weyl} H. Weyl, \textit{The Concept of a Riemann Surface}, 3rd Edition, Addison-Wesley (1955).  Reissued, Dover (2009).

\bibitem{Xia} E.X.W. Xia, ``Congruences Modulo 9 and 27 for Overpartitions," \textit{Ramanujan Journal} 42, pp. 301-323 (2017).

\bibitem{Zuckerman} H.S. Zuckerman, ``Identities Analogous to Ramanujan's Identities Involving the Partition Function," \textit{Duke Mathematical Journal} 5, pp. 88-110 (1939).

\end{thebibliography}
\end{document}